\documentclass[11pt,a4paper]{article}

\usepackage{graphicx,amsfonts,amsbsy, amssymb, amsthm, amsmath}
\usepackage{tikz}
\usetikzlibrary{decorations.pathreplacing}
\usepackage{hyperref,fouridx}
\newcommand{\name}[1]{{#1}}
\newcommand{\shew}{show}

\renewcommand{\vec}[1]{{\bf{#1}}}

\renewcommand{\Re}{\mathop{\mathfrak{Re}}}
\renewcommand{\Im}{\mathop{\mathfrak{Im}}}

\newcommand{\rmd}{{\mathrm d}}
\newcommand{\rme}{{\mathrm e}}
\newcommand{\rmi}{{\mathrm i}}

\newcommand{\Ord}{{\mathrm O}}
\newcommand{\littleo}{{\mathrm o}}

\DeclareSymbolFont{lettersA}{U}{pxmia}{m}{it}

\DeclareMathSymbol{\alphaup}{\mathord}{lettersA}{"0B}
\DeclareMathSymbol{\betaup}{\mathord}{lettersA}{"0C}
\DeclareMathSymbol{\gammaup}{\mathord}{lettersA}{"0D}
\DeclareMathSymbol{\deltaup}{\mathord}{lettersA}{"0E}
\DeclareMathSymbol{\epsilonup}{\mathord}{lettersA}{"22}
\DeclareMathSymbol{\zetaup}{\mathord}{lettersA}{"10}
\DeclareMathSymbol{\etaup}{\mathord}{lettersA}{"11}
\DeclareMathSymbol{\thetaup}{\mathord}{lettersA}{"12}
\DeclareMathSymbol{\iotaup}{\mathord}{lettersA}{"13}
\DeclareMathSymbol{\kappaup}{\mathord}{lettersA}{"14}
\DeclareMathSymbol{\lambdaup}{\mathord}{lettersA}{"15}
\DeclareMathSymbol{\muup}{\mathord}{lettersA}{"16}
\DeclareMathSymbol{\nuup}{\mathord}{lettersA}{"17}
\DeclareMathSymbol{\xiup}{\mathord}{lettersA}{"18}
\DeclareMathSymbol{\piup}{\mathord}{lettersA}{"19}
\DeclareMathSymbol{\rhoup}{\mathord}{lettersA}{"1A}
\DeclareMathSymbol{\sigmaup}{\mathord}{lettersA}{"1B}
\DeclareMathSymbol{\tauup}{\mathord}{lettersA}{"1C}
\DeclareMathSymbol{\upsilonup}{\mathord}{lettersA}{"1D}
\DeclareMathSymbol{\phiup}{\mathord}{lettersA}{"1E}
\DeclareMathSymbol{\chiup}{\mathord}{lettersA}{"1F}
\DeclareMathSymbol{\psiup}{\mathord}{lettersA}{"20}
\DeclareMathSymbol{\omegaup}{\mathord}{lettersA}{"21}

\renewcommand{\Psi}{\varPsi}
\renewcommand{\Lambda}{\varLambda}
\renewcommand{\Sigma}{\varSigma}
\renewcommand{\Gamma}{\varGamma}
\renewcommand{\Theta}{\varTheta}
\renewcommand{\Xi}{\varXi}
\renewcommand{\Pi}{\varPi}
\renewcommand{\Upsilon}{\varUpsilon}
\renewcommand{\Phi}{\varPhi}
\renewcommand{\Omega}{\varOmega}

\newcommand{\R}{{\mathbb R}}
\newcommand{\N}{{\mathbb N}}

\newcommand{\Z}{{\mathbb Z}}

\newcommand{\C}{{\mathbb C}}

\newcommand{\sgn}{\mathop{\rm sign}}

\newcommand{\coloneq}{\mathbin{\hbox{\raise0.08ex\hbox{\rm :}}\!\!=}}
\newcommand{\eqcolon}{\mathbin{=\!\!\hbox{\raise0.08ex\hbox{\rm :}}}}

\renewcommand{\leq}{\leqslant}
\renewcommand{\geq}{\geqslant}
\renewcommand{\epsilon}{\varepsilon} 
\newcommand{\I}{{1\!\!1}}

\newcommand{\dimostrazione}{\noindent{\sl Proof.}\phantom{X}}
\newcommand{\dimostrazionea}[1]{\noindent{\sl Proof of #1.}\phantom{X}}
\newcommand{\finire}{\hspace*{\fill}~$\Box$}

\newcommand \printdate[3]{%
    \def \@suffix##1{%
        \def \@n{##1}%
        \ifnum \@n = 1 st\else%
        \ifnum \@n = 2 nd\else%
        \ifnum \@n = 3 rd\else%
        \ifnum \@n = 21 st\else%
        \ifnum \@n = 22 nd\else%
        \ifnum \@n = 23 rd\else%
        \ifnum \@n = 31 st\else%
        th\fi \fi \fi \fi \fi \fi \fi%
    }%
    \relax%
    \number #1\raise0.7ex\hbox{\footnotesize \@suffix{#1}}\kern0.25em%
    \ifcase #2\or%
        January\or February\or March\or%
        April\or May\or June\or%
        July\or August\or September\or%
        October\or November\or December%
    \fi\ %
    \number #3%
}

\newtheorem{theorem}{Theorem}[section]
\newtheorem{proposition}[theorem]{Proposition}

\newtheorem{lemma}[theorem]{Lemma}

\newcommand{\UN}{{\mathrm U}(N)}
\newcommand{\vand}[1]{\Delta(#1)}
\newcommand{\binomial}[2]{\left(\begin{array}{c} #1 \\ #2\end{array}\right)}
\newcommand{\onefone}[3]{\vphantom{F}_1F_1(#1; #2; #3)}

\newcommand{\oneFone}[1]{\fourIdx{}{1}{(#1)}{1}{F}}
\newcommand{\twoFzero}[1]{\fourIdx{}{2}{(#1)}{0}{F}}
\newcommand{\twoFone}[1]{\fourIdx{}{2}{(#1)}{1}{F}}
\newcommand{\pFq}{\fourIdx{}{p}{(\sigma)}{q}{F}}
\newcommand{\curlyC}{{\mathcal C}}
\newcommand{\curlyW}{{\mathcal W}}
\newcommand{\curlyZ}{{\mathcal Z}}

\numberwithin{equation}{section}

\begin{document}
\title{Derivative moments for characteristic polynomials from the CUE}
\author{B.\ Winn\\{\protect\small\em  Department of Mathematical Sciences, 
Loughborough University,}\\{\protect\small\em Loughborough,
LE11 3TU, U.K. }}
\date{\printdate{1}{9}{2011}
}

\maketitle
\begin{abstract}
  We calculate joint moments of the characteristic polynomial of a
  random unitary matrix from the circular unitary ensemble and its
  derivative in the case that the power in the moments is an odd
  positive integer. The calculations are carried out for finite
  matrix size and in the limit as the size of the matrices goes to
  infinity.  The latter asymptotic calculation allows us to prove a
  long-standing conjecture from random matrix theory.
\end{abstract}
 
\thispagestyle{empty}

\section{Introduction}

There is a deep, and still only partially-understood, relationship 
between analytic number theory, and the theory of random matrices. This
connection goes back to \name{Montgomery} \cite{mon:tpc} who conjectured
that statistical properties of non-trivial zeros of the Riemann
zeta function could be predicted by studying the large $N$ asymptotics
of correlation functions of eigenvalues of $N\times N$ random unitary
matrices. This conjecture is supported by theoretical \cite{mon:tpc,
hej:ott, rud:tnl, rud:zop}, heuristic \cite{bog:rmtI,bog:rmtII, bog:gtf} 
and numerical \cite{odl:otd, odl:zeta, hia:tzf} evidence.

The full power of the conjectured relationships between random matrix
theory and number theory is found in the study of moments of the
Riemann zeta function.  Using random matrix theory, mathematicians
have been able to make predictions for moments of various kinds, where
no conjectures or guesses existed before \cite{kea:rmtz, kea:rmtL,
  hug:rmt, hug:otc, con:aor, con:imo, con:lot}.

A number of review articles have appeared, such as \cite{con:Lfa,
  kea:rmaL, kea:rma, sna:rza}, which summarise the main developments
that have occurred over the past few years.

The main object of our study will be the $N$-dimensional circular
unitary ensemble (CUE) of random matrix theory.  This is the
probability space consisting of the set $\UN$ of $N\times N$ unitary
matrices, equipped with normalised Haar measure, $\mu^{\rm Haar}$.

For a matrix $U\in \UN$ we denote the characteristic polynomial by
\begin{equation} \label{eq:def_Z}
  Z_U(\theta) \coloneq \prod_{n=1}^N \left( 1-\rme^{\rmi(\theta_n-\theta)}
\right),
\end{equation}
where $\rme^{\rmi\theta_1},\ldots,\rme^{\rmi\theta_N}$ are the eigenvalues
of $U$.

Define
\begin{equation} \label{eq:def_V}
  V_U(\theta)\coloneq \exp\!\left(\rmi N\frac{\theta+\pi}2 -\rmi
{\sum}_{n=1}^N \frac{\theta_n}2\right) Z_U(\theta).
\end{equation}
Then $V_U(\theta)$ is real-valued for $\theta\in[0,2\pi)$.

In recent years there has been interest in the joint moments of the
distribution of $V_U$ and its derivative. Define, for $h>-1/2$ and
$k>h-1/2$,
\begin{equation} \label{eq:def_Ftwiddle}
  \tilde{F}_N(h,k)\coloneq \int_{\UN} |V_U(0)|^{2k-2h} |V'_U(0)|^{2h}\,
\rmd\mu^{\rm Haar},
\end{equation}
and the limiting values
\begin{equation} \label{eq:moments}
\tilde{F}(h,k)\coloneq \lim_{N\to\infty} \frac1{N^{k^2+2h}}
\tilde{F}_N (h,k).
\end{equation}
When $h=0$, the moments \eqref{eq:def_Ftwiddle} of $V_U$ are
precisely the same as the moments of the characteristic polynomial.
\name{Keating} and \name{Snaith} \cite{kea:rmtz} considered
$\tilde{F}_N(0,k)$, and proved that
\begin{equation}
  \label{eq:keating_snaith}
  \tilde{F}_N(0,k) =\prod_{j=1}^N \frac{\Gamma(j)\Gamma(j+2k)}
{\Gamma(j+k)^2},
\end{equation}
and \shew ed that \eqref{eq:keating_snaith} extends to
the region $\Re\{k\} > -1/2$ of the complex plane.

Let $\curlyZ(t)$ denote Hardy's function:
\begin{equation}
 \curlyZ(t)\coloneq \rme^{\rmi\vartheta(t)}\zeta({\textstyle\frac12}+
\rmi t),
\end{equation}
where
\begin{equation}
 \vartheta(t)\coloneq \Im\left\{\log\left( \pi^{-\rmi t/2}\Gamma\left(
      {\textstyle\frac14 + \frac12\rmi t}\right)\right)\right\},
\end{equation}
and $\zeta(s)$ and $\Gamma(s)$ denote respectively the Riemann zeta
function and the Euler gamma function. It follows from the functional
equation for $\zeta(s)$ that $\curlyZ(t)$ is real for $t\in\R$, and
it is apparent that
 $|\zeta({\textstyle\frac12}+\rmi t)| = |\curlyZ(t)|$, so
$\curlyZ(t)$ is to the Riemann zeta function as $V_U(\theta)$ is to
the characteristic polynomial of a random unitary matrix. A series
of conjectures due to \name{Hall} \cite{hal:awt},
\name{Conrey} and \name{Ghosh} \cite{con:omvI},
and \name{Hughes} \cite{hug:phd} has culminated in the following prediction
for joint moments of $\curlyZ(t)$ and its derivative:
\begin{equation}
  \label{eq:rmt_conjecture}
  \frac1T \int_0^T |\curlyZ(t)|^{2k-2h}|\curlyZ'(t)|^{2h}\,\rmd t
  \sim\tilde{F}(h,k) A(k) \left( \log T\right)^{k^2+2h},
\qquad\mbox{as $T\to\infty$,}
\end{equation}
where
\begin{equation}
  \label{eq:a}
  A(k)\coloneq \prod_{p\text{ prime}} \left(1-\frac1p\right)^{k^2}
\sum_{\ell=0}^\infty
\left( \frac{\Gamma(\ell+k)}{\ell!\Gamma(k)}\right)^2 p^{-\ell}.
\end{equation}

\name{Hughes} \cite{hug:phd} used random matrix theory to calculate
$\tilde{F}(h,k)$ for $h=1,2,3$ and \name{Dehaye} \cite{deh:jmo,
  deh:ano} has derived formul\ae\ for $\tilde{F}(h,k)$ for all
$h\in\N$ in terms of sums over partitions (see section
\ref{sec:main-results} below for a precise statement). Using their
results, the values $\tilde{F}(1,1)=1/12$, $\tilde{F}(1,2)=1/720$ and
$\tilde{F}(2,2)=1/6720$ can be calculated. The corresponding moments
\eqref{eq:rmt_conjecture} for Hardy's function have been calculated by
\name{Ingham} \cite{ing:mvt} and \name{Conrey} \cite{con:tfm}, and
give complete agreement for these values of $k$ and $h$.
\nopagebreak

\name{Conrey} and \name{Ghosh} \cite{con:omvII} have also proved 
(assuming the Riemann hypothesis) that
\begin{equation}
  \label{eq:conrey_ghosh}
  \frac1T\int_1^T |\curlyZ(t)\curlyZ'(t)|\,\rmd t \sim
\frac{\rme^2-5}{4\pi}\left(\log T\right)^2.
\end{equation}
This is proved by relating the joint moment to a discrete second
moment of the Riemann zeta function at its successive extrema on the
critical line, which had been calculated earlier in
\cite{con:amv}. The latter result was proved by an integration against
the logarithmic derivative of a function with zeros at the locations
of maxima of $|\zeta(\frac12+\rmi t)|$ and which could be
well-approximated by a Dirichlet series. The numerical constants in
\eqref{eq:conrey_ghosh} arise as values of residues at poles in the
relevant contour integral.

The asymptotic \eqref{eq:conrey_ghosh} naturally leads to the 
conjecture \cite[page 110]{hug:phd} that
\begin{equation}
  \label{eq:conjecture}
  \tilde{F}\left( {\textstyle\frac12}, 1\right) = \frac{\rme^2-5}{4\pi}.
\end{equation}
However, most attention on the problem of calculating moments 
\eqref{eq:def_Ftwiddle} has focussed on integer values of $h$. In this work
we will take the first steps beyond integer values of $h$, by studying
$\tilde{F}_N(h,k)$ for half-integer values of $h$. In particular we will
supply a proof of \eqref{eq:conjecture}.

\section{Main results} \label{sec:main-results}

In order to put our results into context, we first recall a
result of \name{Dehaye} \cite{deh:jmo, deh:ano}. To do this it will be 
necessary to fix some notations regarding combinatorics of partitions.

We recall that a {\em partition\/} is a finite sequence $\lambda =
(\lambda_1,\ldots,\lambda_j)$ with $\lambda_1\geq \lambda_2\geq \cdots
\geq \lambda_j$. $j$ is the number of parts of $\lambda$, which is
also denoted by $\ell(\lambda)$.
The sum of the parts of $\lambda$ is denoted by
$|\lambda|=\lambda_1+\cdots+\lambda_j$. For $n\in\N$ we write
$\lambda \vdash n$ if $|\lambda|=n$, 
and $\lambda \vdash_k n$ if $|\lambda|=n$ with $\ell(\lambda)\leq k$:
a partition
of $n$ into not more than $k$ parts.

The {\em generalised Pochhammer symbol} $[b]_\lambda^{(\sigma)}$ is 
defined for a partition $\lambda$, a parameter $\sigma>0$ and $b\in\R$ by
\begin{equation}
  \label{eq:1}
  [b]_\lambda^{(\sigma)} \coloneq
\prod_{i=1}^{\ell(\lambda)}\prod_{j=1}^{\lambda_i}
\left(b + j -1 - \frac{i-1}\sigma\right).
\end{equation}
We will most often be taking the parameter $\sigma=1$, so that we
define the special notation
\begin{equation}
  \label{eq:29}
  [b]_\lambda \coloneq [b]_\lambda^{(1)}.
\end{equation}
In terms of the usual (rising) Pochhammer symbol $(\cdot)_\cdot$, we have
\begin{equation}
  [b]_\lambda = \prod_{i=1}^{\ell(\lambda)} (b-i+1)_{\lambda_i}.
\end{equation}

\begin{figure}
  \centering
  \begin{tikzpicture}
    \foreach \position in {(0,0),(0,0.5),(0,1),(0,1.5),(0.5,1),(0.5,1.5),
      (1,1),(1,1.5),(1.5,1.5)}
    \draw \position rectangle +(0.5,0.5);
    \draw (4,0) rectangle +(0.5,0.5) node [anchor=north east] {$1$};
    \draw (4,0.5) rectangle +(0.5,0.5) node [anchor=north east] {$2$};
    \draw (4,1) rectangle +(0.5,0.5) node [anchor=north east] {$5$};
    \draw (4,1.5) rectangle +(0.5,0.5) node [anchor=north east] {$7$};
    \draw (4.5,1) rectangle +(0.5,0.5) node [anchor=north east] {$2$};
    \draw (4.5,1.5) rectangle +(0.5,0.5) node [anchor=north east] {$4$};
    \draw (5,1) rectangle +(0.5,0.5) node [anchor=north east] {$1$};
    \draw (5,1.5) rectangle +(0.5,0.5) node [anchor=north east] {$3$};
    \draw (5.5,1.5) rectangle +(0.5,0.5) node [anchor=north east] {$1$};
  \end{tikzpicture}
  \caption{The Ferrers diagram (left) and the hook-lengths (right) for the partition $\lambda=(4,3,1,1)$.}
  \label{fig:ferrers}
\end{figure}
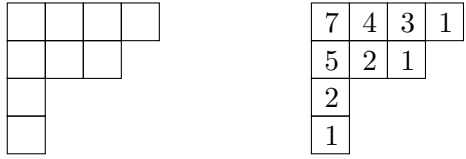

A partition can be represented graphically by a {\em Ferrers
  diagram\/} (see figure \ref{fig:ferrers}), in which parts of a
partition are represented by a vertical arrangement of boxes aligned
at the left-hand side. For each box $\Box$ in the Ferrers diagram, the
{\em arm-length\/} $a(\Box)$ is the number of boxes strictly to the
right of $\Box$, and the {\em leg-length} $g(\Box)$ is the number of
boxes strictly below $\Box$. The {\em hook-length} of $\Box$ is
$a(\Box)+g(\Box)+1$: the number of boxes to the right and below, with
the box itself counted exactly once.  For the partition
$\lambda=(4,3,1,1)$, the hook-lengths are indicated in figure
\ref{fig:ferrers}. The product of all hook-lengths will be called the
hook-length of the partition, and denoted $h_\lambda$.  For example,
for $\lambda=(4,3,1,1)$ we find $h_\lambda=1680$.

The Ferrers diagram for a partition can be used to define the
{\em transpose} partition, by reflection of the diagram about the
main diagonal. So for $\lambda=(4,3,1,1)$, the transpose
partition is $\lambda^{\rm T}=(4,2,2,1)$. Clearly the length
$\ell(\lambda^{\rm T})$ of a transpose partition
is equal to the size of the largest part of $\lambda$,
and $|\lambda^{\rm T}| = |\lambda|$. 
It is also straightforward to see that
\begin{equation}
  \label{eq:106}
  [b]_{\lambda^{\rm T}} = (-1)^{|\lambda|}[-b]_\lambda
\qquad\mbox{and}\qquad 
  h_{\lambda^{\rm T}}=h_\lambda.
\end{equation}

We define two quantities in terms of a sum over partitions. Let
$k\in\N$ and $p\in\N_0\coloneq \N\cup\{0\}$. Then 
\begin{equation}
 \label{eq:curlyCN}
  \curlyC_N(p,k) \coloneq (-2)^p \sum_{\lambda\vdash_k p} \frac{[k]_\lambda 
[-N]_{\lambda}}{[2k]_\lambda h_\lambda^2},
\end{equation}
and
\begin{equation}
  \label{eq:curlyC}
 \curlyC(p,k) \coloneq 2^p \sum_{\lambda\vdash_k p} \frac{[k]_{\lambda}}
{[2k]_\lambda h_{\lambda}^2}.
\end{equation}
We observe that 
\begin{equation}\label{eq:curlyC_asympt}
\curlyC_N(p,k)\sim \curlyC(p,k) N^p\qquad\mbox{as $N\to\infty$.}
\end{equation}
Related quantities appear in the work of \name{Dehaye} (see the comment
following theorem \ref{thm:dehaye} below).

The main result of \name{Dehaye} \cite{deh:jmo, deh:ano} relevant to
our work is the following:
\begin{theorem} \label{thm:dehaye}
  For $h,k\in\N$ with $k>h-\frac12$,
  \begin{equation}
    \label{eq:2}
    \tilde{F}_N(h,k) = \frac{(-1)^h}{2^{2h}} \tilde{F}_N(0,k)
    \sum_{p=0}^{2h} \frac{(2h)!(-N)^{2h-p}}{(2h-p)!} \curlyC_N(p,k).
  \end{equation}
\end{theorem}
Moreover, for fixed $h\in\N$, \name{Dehaye} has \shew n that
the equation \eqref{eq:2} extends meromorphically
in $k$ to the region $\Re\{k\}>h-\frac12$ of the complex plane.

We have re-written \name{Dehaye}'s result using our notation. In fact
he considers quantities related to $\curlyC_N$ and $\curlyC$ (see for
example equations (10) and (11) of \cite{deh:ano}), which are defined
similarly to \eqref{eq:curlyCN} and \eqref{eq:curlyC}, but without the
restriction on the number of parts of $\lambda$ in the summation. The
presence of the factor $[k]_\lambda$ mean that his and our quantities
coincide for $p\leq 2k$ (but could be different for $p>2k$). Since the
sum in \eqref{eq:2} goes only up to $2h$, this difference is not
pertinent in theorem \ref{thm:dehaye}.

The main result of our work is the following, which gives an explicit
formula for $\tilde{F}_N(h,k)$ for half-integer $h$:
\begin{theorem} \label{thm:main}
  Let $h=(2m-1)/2$ for $m\in\N$ and let $k\in\N$ with $k>h-\frac12$.
Then
\begin{multline}
  \label{eq:3}
  \tilde{F}_N(h,k) =\frac{2(-1)^{h+1/2}}{2^{2h}\pi} \tilde{F}_N(0,k)
\Bigg\{ \sum_{p=1}^{2h}\sum_{\ell=1}^p
 \binomial{2h}{p-\ell} \frac{(-1)^\ell}\ell
(-N)^{2h-p} p! \curlyC_N(p,k) \\
+\sum_{p=2h+1}^{kN} \frac{(2h)!(p-2h-1)!}{N^{p-2h}} \curlyC_N(p,k)\Bigg\}.
\end{multline}
\end{theorem}

The paper is structured as follows: In section \ref{sec:drei} we
write down an integral representation for $\tilde{F}_N(h,k)$ involving
an integration over a real parameter $\zeta$ and a multi-dimensional
integral of size $N$. In section \ref{sec:vier} we evaluate the
multi-dimensional integral in closed form, and in section
\ref{sec:some-integr-involv} we calculate asymptotics of some integrals
related to the $\zeta$-integral. We give the proof of theorem
\ref{thm:main} in section \ref{sec:sechs}, and we use the theorem to give
evaluations of some moments (including \eqref{eq:conjecture}) in
section \ref{sec:partition-sums-proof}. 

\section{An integral representation}\label{sec:drei}

\subsection{Notation and properties of Vandermonde determinants}
Let us fix some notation which will be used throughout the remainder of
the paper. Let $\vec{x}=(x_1,\ldots,x_N)$. We will denote in multiple
integrals,
\begin{equation}
  \label{eq:dvecx}
 \rmd\vec{x} = \rmd x_1\cdots \rmd x_N.
\end{equation}

We shall also denote by
\begin{equation}
  \label{eq:vandermonde}
  \vand{\vec{x}}\coloneq  \prod_{1\leq j<k\leq N} (x_k-x_j),
\end{equation}
the Vandermonde determinant. It seems prudent at this stage to note a
few properties of $\vand{\cdot}$, that we will make use of later
\cite{ait:dam}.

First of all, note that the Vandermonde determinant is a matrix determinant.
We have
\begin{align}
  \label{eq:vandermonde:2}
  \vand{\vec{x}} &= \det\left( x_i^{j-1}\right)_{i,j=1,\ldots,N}\\
&=\sum_{\sigma\in S_N}\sgn(\sigma) x_1^{\sigma(1)-1}\cdots 
x_N^{\sigma(N)-1},
\label{eq:vandermonde:4}
\end{align}
where $S_N$ is the symmetric group on $N$ elements.
We see, therefore, that $\vand{\vec{x}}$ is a homogeneous polynomial
in the variables $x_1,\ldots, x_N$, of degree $N(N-1)/2$.

If $\{p_j(\cdot)\}_{j=1,\ldots,N}$ is a set of monic polynomials with
the degree of $p_j$ being $j-1$, then an alternative expression for
$\vand{\vec{x}}$ is
\begin{equation}
  \label{eq:vandermonde:3}
  \vand{\vec{x}}=\det\left( p_j(x_i) \right)_{i,j=1,\ldots,N},
\end{equation}
which may be proved by applying elementary column operations to the 
representation \eqref{eq:vandermonde:2}.

In order to justify the convergence of certain integrals, we shall employ
the following crude bound.
\begin{lemma}
  \label{lem:vandermonde}
For $\vec{x}\in\R^N$,
\begin{equation}
  |\vand{\vec{x}}|\leq N! \left( (1+x_1^2)^{1/2}\cdots
  (1+x_N^2)^{1/2}\right)^{N-1}.
\end{equation}
\end{lemma}
\dimostrazione
For a permutation $\sigma$,  $\sgn\sigma=\pm1$. So,
by \eqref{eq:vandermonde:4} we can bound
\begin{align}
  |\vand{\vec{x}}|&\leq \sum_{\sigma\in S_N} |x_1|^{\sigma(1)-1}\cdots
|x_N|^{\sigma(N)-1} \nonumber \\
&\leq \Big(\sum_{\sigma\in S_N} 1\Big) (1+x^2_1)^{(N-1)/2}\cdots
(1+x_N^2)^{(N-1)/2},
\end{align}
and use the fact that $S_N$ has order $N!$. 
\finire

\mathversion{bold}
\subsection{An integral representation for 
\texorpdfstring{$\tilde{F}_N(h,k)$}{FtwiddleN(h,k)}}
\mathversion{normal}

To evaluate the averages $\tilde{F}_N(h,k)$ we express this quantity as a
multi-dimensional integral:
\begin{proposition} \label{prop:zwei}
  Let $n\in\N_0$, and define
\begin{equation} \label{eq:def_K}
K_n(\epsilon,\zeta)\coloneq \frac{(-1)^n}{\pi}\frac{\partial^n}{\partial 
\epsilon^n}
\left(\frac{\epsilon}{\epsilon^2+\zeta^2}\right).
\end{equation}
Then if $2h\in\N_0$ and $k>h-\frac12$,
\begin{equation}
  \label{eq:moment_multi}
  \tilde{F}_N(h,k) = \lim_{\epsilon\downarrow 0}\frac{2^{N^2+2kN-2h}}{(2\pi)^N
N!}\int_{-\infty}^\infty  \int_{-\infty}^{\infty} \!\cdots\!
\int_{-\infty}^{\infty} K_{2h}(\epsilon,\zeta)
\prod_{j=1}^N\frac{\rme^{\rmi \zeta x_j}}
{(1+x_j^2)^{N+k}} \vand{\vec{x}}^2\,\rmd\vec{x}
\rmd\zeta.
\end{equation}
\end{proposition}

In order to prove proposition \ref{prop:zwei}, let us collect a few
auxiliary results.

\begin{lemma} \label{lem:due}
  For $x\in\R$ and $\epsilon>0$, we have
\begin{equation}
\int_{-\infty}^\infty K_n(\epsilon,\zeta)\rme^{\rmi x \zeta}\,\rmd\zeta
=|x|^n\rme^{-\epsilon|x|},
\end{equation}
and the integral converges uniformly in $x$ and $\epsilon\geq \epsilon_0
>0$.
\end{lemma}
\dimostrazione
By calculating residues, or the Fourier inversion theorem, we know that
for $\epsilon>0$,
\begin{equation}
  \frac1\pi\int_{-\infty}^\infty \frac{\epsilon}{\epsilon^2+\zeta^2}
\rme^{\rmi x \zeta}\,\rmd\zeta = \rme^{-\epsilon|x|}, 
\end{equation}
uniformly for $\epsilon\geq \epsilon_0>0$.
To justify differentiation under the integral, we note that
\begin{equation} \label{eq:K_estimate}
  K_n(\epsilon,\zeta)\ll \frac1{(\epsilon^2+\zeta^2)^{(n+1)/2}}
\end{equation}
for $n\geq1$, so that uniform convergence of the resulting integrals is
assured.
\finire

To compute averages over $\UN$, the most useful tool available is
Weyl's integration formula.
For any function $f(U)$ of a unitary matrix,
which depends only on the $N$ eigenvalues $\rme^{\rmi\theta_1},\ldots,
\rme^{\rmi \theta_N}$,
\begin{equation}
  \label{eq:weyl_integration}
  \int_{\UN}f(U)\,\rmd\mu^{\rm Haar} = \frac1{(2\pi)^N N!}\int_0^{2\pi}
\!\cdots\int_0^{2\pi} f(U)\prod_{1\leq j < k \leq N}\left| \rme^{\rmi\theta_k}
-\rme^{\rmi\theta_j}\right|^2\,\rmd\theta_1\cdots\rmd\theta_N.
\end{equation}
By following the substitutions made in \cite{kea:rmtz}, we can express
$\tilde{F}(h,k)$ as an integral over $\R^N$.

\begin{proposition}
  It follows from \eqref{eq:weyl_integration} that
\begin{equation}\label{eq:zero}
\tilde{F}_N(h,k)=\frac{2^{N^2+2kN-2h}}{(2\pi)^N N!} 
\int_{-\infty}^{\infty} \!\cdots\!
\int_{-\infty}^{\infty} 
\prod_{j=1}^N\frac{1}{(1+x_j^2)^{N+k}} \left|
x_1+\cdots+x_N\right|^{2h}\vand{\vec{x}}^2\,\rmd\vec{x}.
\end{equation}
\end{proposition}
\dimostrazione
  Differentiating \eqref{eq:def_V}, we get
  \begin{equation}
 V'_U(\theta)=\frac{\rmi N}2 V_U(\theta) + \frac{V_U(\theta)}{Z_U(\theta)}
 Z'_U(\theta).
  \end{equation} 
So,
\begin{equation} \label{eq:uno}
  V'_U(0)=V_U(0)\left(\frac{\rmi N}2 + \frac{Z'_U(0)}{Z_U(0)}\right).
\end{equation}
Furthermore,
\begin{align}
\frac{Z_U'(\theta)}{Z_U(\theta)}&=\frac{\rmd}{\rmd\theta}\log(Z_U(\theta))
\nonumber \\
&=\sum_{n=1}^N\frac{\rmi\rme^{\rmi(\theta_n-\theta)}}
{1-\rme^{\rmi(\theta_n-\theta)}},
\end{align}
from \eqref{eq:def_Z}. We then get
\begin{align}
  V'_U(0)&=\rmi V_U(0)\Big(\frac N2 + \sum_{n=1}^N \frac1{\rme^{-\rmi\theta_n}
-1}\Big) \nonumber\\ \label{eq:due}
&=-\frac12 V_U(0)\sum_{n=1}^N\cot\!\left(\frac{\theta_n}2\right),
\end{align}
using the fact that
\begin{equation}
\frac1{\rme^{-\rmi\theta_n}-1} = \frac{\rmi}2\cot\!\left(\frac{\theta_n}2
\right)-\frac12.
\end{equation}

Substituting \eqref{eq:due} into \eqref{eq:def_Ftwiddle}, we get
\begin{equation}
  \tilde{F}(h,k) =\frac1{2^{2h}}\int_{\UN} |V_U(0)|^{2k}\Bigg|\sum_{n=1}^N
\cot\!\left(\frac{\theta_n}2\right)\Bigg|^{2h}\,\rmd\mu^{\rm Haar}.
\end{equation}
Following \cite{kea:rmtz}, we write the integral over $\UN$ as a
multiple-integral using Weyl's identity \eqref{eq:weyl_integration},
and make the substitutions $x_j=\cot(\theta_j/2)$ therein, leading to
\eqref{eq:zero}.  \finire

\dimostrazionea{proposition \ref{prop:zwei}}
By lemma \ref{lem:due}, for any $\epsilon>0$,
\begin{multline} \label{eq:tre}
 \int_{-\infty}^{\infty} \!\cdots\!
\int_{-\infty}^{\infty} \int_{-\infty}^\infty K_{2h}(\epsilon,\zeta)
\prod_{j=1}^N\frac{\rme^{\rmi \zeta x_j}}
{(1+x_j^2)^{N+k}} \vand{\vec{x}}^2\,\rmd\zeta\rmd\vec{x}
=\\
 \int_{-\infty}^{\infty} \!\cdots\! \int_{-\infty}^{\infty}
\prod_{j=1}^N\frac{1}
{(1+x_j^2)^{N+k}} \left|x_1+\cdots+x_N\right|^{2h}
\rme^{-\epsilon|x_1+\cdots+x_N|} \vand{\vec{x}}^2\,\rmd\vec{x}.
\end{multline}
By lemma \ref{lem:vandermonde} and equation \eqref{eq:K_estimate} we
see that the left-hand side of \eqref{eq:tre} is absolutely
integrable, so by the Fubini-Tonelli theorem, we can move the
$\zeta$-integral to the outside.

Considering the right-hand side of \eqref{eq:tre}, we estimate
\begin{align}
  |x_1+\cdots+x_N|^{2h}&\leq \left( \sqrt{1+x_1^2}+\cdots+\sqrt{1+x_N^2}
\right)^{2h} \nonumber \\
&\leq N^{2h}\left( (1+x_1^2)^{1/2}\cdots(1+x_N^2)^{1/2}\right)^{2h}
\nonumber \\
&=N^{2h}(1+x_1^2)^h\cdots(1+x^2_N)^h.
\end{align}
Combining this with lemma \ref{lem:vandermonde}, we get
\begin{equation}
\prod_{j=1}^N\frac{1}
{(1+x_j^2)^{N+k}} \left|x_1+\cdots+x_N\right|^{2h}
 \vand{\vec{x}}^2 \leq (N^h N!)^2\prod_{j=1}^N\frac1{(1+x_j^2)^{k+1-h}},
\end{equation}
so that if $k>h-1/2$, the right-hand side of \eqref{eq:tre} is uniformly
convergent in $\epsilon$, and we may pass the limit $\epsilon\downarrow 0$
under the integration. The resulting equality is \eqref{eq:moment_multi}.
\finire

\section{Evaluation of the multi-dimensional integrals}\label{sec:vier}

The main calculation in this paper is an evaluation of the
integral
\begin{equation}
  \label{eq:4}
  H(k,\zeta) \coloneq \int_{-\infty}^\infty \!\cdots\!
\int_{-\infty}^\infty \prod_{j=1}^N \frac{\rme^{\rmi\zeta x_j}}
{(1+x_j^2)^{N+k}}\vand{\vec{x}}^2\,\rmd\vec{x}.
\end{equation}
To this end, we will first derive an equivalent representation for
$H(k,\zeta)$ which we will be able to evaluate, for integer $k$, in
terms of Laguerre polynomials, and multivariable hypergeometric
functions.

\begin{proposition}\label{prop:integral}
  Let $\zeta>0$ be fixed, and $k\in\C$ with $\Re\{ k\}>-\frac12$. Then
  \begin{equation}
    \label{eq:5}
    H(k,\zeta) = \frac{\pi^N}{2^{(N+2k-1)N}}\prod_{j=0}^{N-1}
 \frac1{\Gamma(k+1+j)^2}\rme^{-N\zeta}\int_0^\infty\!\cdots\!\int_0^\infty
\prod_{j=1}^N (y_j+2\zeta)^k y_j^k\rme^{-y_j} \vand{\vec{y}}^2\,
\rmd\vec{y}.
  \end{equation}
\end{proposition}

\dimostrazione 
Let 
\begin{equation}
  \label{eq:6}
  L(k,\zeta)\coloneq
\int_0^\infty\!\cdots\!\int_0^\infty
\prod_{j=1}^N (y_j+2\zeta)^k y_j^k\rme^{-y_j} \vand{\vec{y}}^2\,
\rmd\vec{y}.
\end{equation}

By \eqref{eq:vandermonde:3} we may write, for fixed $\zeta$,
\begin{align}
  \nonumber
\vand{\vec{y}}&=\det\left( (y_i+2\zeta)^{j-1} \right)_{i,j=1,\ldots,N} \\
&=\sum_{\sigma\in S_N}\sgn(\sigma)(y_1+2\zeta)^{\sigma(1)-1}
\cdots(y_N+2\zeta)^{\sigma(N)-1}.
\end{align}
Using the expression \eqref{eq:vandermonde:2} for the second of the
two Vandermonde factors, we can express $L(k,\zeta)$ as a sum of
products of integrals. We get
\begin{equation}
  \label{eq:laguerre_due}
L(k,\zeta) = \sum_{\sigma,\tau \in S_N}\sgn(\sigma)\sgn(\tau)
\prod_{j=1}^N {\mathcal I}_{\sigma(j)-1,\tau(j)-1},
\end{equation}
where
\begin{equation}
  \label{eq:curly_I}
  {\mathcal I}_{\mu,\nu} \coloneq 
\int_0^\infty (y+2\zeta)^{\nu+k}y^{\mu+k}\rme^{-y}\,\rmd y.
\end{equation}

Let us now consider the integral 
\begin{equation}\label{eq:laguerre_quattro}
\int_{-\infty}^{\infty} \!\cdots\! \int_{-\infty}^{\infty} 
\prod_{j=1}^N\frac{\rme^{\rmi \zeta x_j+\zeta}}
{(1+x_j^2)^{N+k}} \vand{\vec{x}}^2\,\rmd\vec{x}.
\end{equation}
We exploit the homogeneity of $\vand{\vec{x}}$ to write
\begin{align} \nonumber
  \prod_{j=1}^N\frac{1}
{(1+x_j^2)^{N-1}} \vand{\vec{x}}^2 &= \prod_{1\leq j<m\leq N} 
\frac{(x_m-x_j)^2}{(1+x_j^2)(1+x_m^2)} \\  \nonumber
&= \prod_{1\leq j<m\leq N}\frac{x_m-x_j}{(1+\rmi x_j)(1+\rmi x_m)}
\prod_{1\leq j<m\leq N}\frac{x_m-x_j}{(1-\rmi x_j)(1-\rmi x_m)}\\ \nonumber
&= \prod_{1\leq j<m\leq N}\left(\frac1{1+\rmi x_m} -\frac{1}{1+\rmi x_j}
\right)
\prod_{1\leq j<m\leq N}\left(\frac1{1-\rmi x_m} -\frac{1}{1-\rmi x_j}
\right)\\
=\sum_{\sigma,\tau\in S_N}&\frac{\sgn(\sigma)\sgn(\tau)}
{(1+\rmi x_1)^{\tau(1)-1}(1-\rmi x_1)^{\sigma(1)-1}\cdots
(1+\rmi x_N)^{\tau(N)-1}(1-\rmi x_N)^{\sigma(N)-1}}.
\end{align}
Substituting this into \eqref{eq:laguerre_quattro} we get
\begin{equation}
  \label{eq:laguerre_cinque}
  \int_{-\infty}^{\infty} \!\cdots\! \int_{-\infty}^{\infty} 
\prod_{j=1}^N\frac{\rme^{\rmi \zeta x_j+\zeta}}
{(1+x_j^2)^{N+k}} \vand{\vec{x}}^2\,\rmd\vec{x} 
=\sum_{\sigma,\tau \in S_N}\sgn(\sigma)\sgn(\tau)
\prod_{j=1}^N {\mathcal J}_{\sigma(j)-1,\tau(j)-1},
\end{equation}
where
\begin{equation} \label{eq:J_integral}
  {\mathcal J}_{\mu,\nu}\coloneq \int_{-\infty}^\infty 
\frac{\rme^{\rmi\zeta x
+\zeta}}{(1+\rmi x)^{1+k+\nu}(1-\rmi x)^{1+k+\mu}} \,\rmd x.
\end{equation}
The integral \eqref{eq:J_integral} can be evaluated in terms of a
confluent hypergeometric\footnote{For the definition of
  $\vphantom{F}_1F_1$, we refer the reader to section
  \ref{sec:eval-terms-hyperg}.} function---from \cite[formula
3.384.9]{gra:tis} we find that for $\alpha,\beta\in\N$ and for
$\zeta>0$,
\begin{equation*}
  \int_{-\infty}^\infty \frac{\rme^{\rmi\zeta x}}{(1+\rmi x)^\alpha
(1-\rmi x)^\beta}\,\rmd x = \frac{\pi}{2^{\alpha+\beta-2}}
\frac{\Gamma(\alpha+\beta-1)}{\Gamma(\alpha)\Gamma(\beta)}\rme^{-\zeta}
\onefone{1-\alpha}{2-\alpha-\beta}{2\zeta}
\end{equation*}
---however, our main concern is to \shew\ that ${\mathcal J}_{\mu,\nu}$ is
equal to ${\mathcal I}_{\mu,\nu}$, up to a constant, which we
do next.

We, temporarily, assume that $\Re\{ k\}>0$. Since
\begin{equation}
  \label{eq:7}
  \frac1{(1-\rmi x)^{1+k+\mu}} = \frac1{\Gamma(1+k+\mu)}\int_0^\infty
\rme^{-(1-\rmi x)u}u^{k+\mu}\,\rmd u,
\end{equation}
we may write
\begin{equation}
  \label{eq:8}
  {\mathcal J}_{\mu,\nu} = \frac1{\Gamma(1+k+\mu)}\int_{-\infty}^\infty
\int_0^\infty \rme^{-(1-\rmi x)u}u^{k+\mu} \frac{\rme^{\rmi\zeta x-\zeta}}
{(1+\rmi x)^{1+k+\nu}}\,\rmd u\rmd x.
\end{equation}
We bound
\begin{equation}
  \label{eq:9}
  \left| \frac{1}{(1+\rmi x)^{1+k+\nu}}\right| \leq \frac{\rme^{(\pi/2)
\Im\{k\}}}{(1+x^2)^{(1+\nu+\Re\{k\})/2}},
\end{equation}
and
\begin{equation}
  \label{eq:10}
  \left| \rme^{-(1-\rmi x)u}u^{k+\mu}\right| \leq u^{\Re\{k\}+\mu}
\rme^{-u},
\end{equation}
so that the double integral \eqref{eq:8} is absolutely convergent and
we may reverse the order of integration. Therefore, we have
\begin{align}
  \label{eq:11}
  {\mathcal J}_{\mu,\nu} &= \frac1{\Gamma(1+k+\mu)}\int_0^\infty
  \rme^{\zeta-u}u^{k+\mu}\int_{-\infty}^\infty \frac{\rme^{\rmi(\zeta
+u)x}}{(1+\rmi x)^{1+k+\nu}}\,\rmd x\,\rmd u \nonumber \\
&=\frac{2\pi}{\Gamma(1+k+\mu)\Gamma(1+k+\nu)}\int_0^\infty
\rme^{\zeta-u}u^{k+\mu}(u+\zeta)^{k+\nu}\rme^{-(\zeta+u)}\,\rmd u
\nonumber \\
&=\frac{\pi}{2^{2k+\mu+\nu}\Gamma(1+k+\mu)\Gamma(1+k+\nu)}\int_0^\infty
(y+2\zeta)^{k+\nu}y^{k+\mu}\rme^{-y}\,\rmd y
\nonumber \\
&=\frac{\pi}{2^{2k+\mu+\nu}\Gamma(1+k+\mu)\Gamma(1+k+\nu)}
{\mathcal I}_{\mu,\nu},
\end{align}
using Laplace's formula (equation 3.382.6 of \cite{gra:tis})
\begin{equation}
  \label{eq:laplace}
  \frac1{2\pi}\int_{-\infty}^\infty \frac{\rme^{\rmi p x}}{(1+\rmi x)^s}
\,\rmd x = \frac{p^{s-1} \rme^{-p}}{\Gamma(s)}, \qquad
\mbox{$p>0$,}
\end{equation}
to pass from the first to the second line of \eqref{eq:11}.

The relationship \eqref{eq:11}, together with \eqref{eq:laguerre_cinque}
and \eqref{eq:laguerre_due} proves \eqref{eq:5} for $\Re\{ k\}>0$. To
complete the proof we note that both sides of \eqref{eq:5} may be
continued as analytic functions of $k$ to $\Re\{k\}>-\frac12$. \finire

For integer values of $k$ we are able to give two direct evaluations
of the integral in the right hand side of \eqref{eq:5}. The
first one uses Laguerre polynomials, and the second uses 
a hypergeometric function of matrix argument.

\subsection{Evaluation in terms of Laguerre polynomials}
We recall that the classical Laguerre polynomials $L_n^{(\alpha)}$
are defined for a parameter $\alpha>-1$ by
\begin{equation}
  \label{eq:laguerre}
  L_n^{(\alpha)}(t) \coloneq  \frac{\rme^t}{t^\alpha n!}
\frac{\rmd^n}{\rmd t^n}\left(t^{\alpha+n}\rme^{-t}\right).
\end{equation}
An explicit formula for $L_n^{(\alpha)}$ is given by
\begin{equation}
  \label{eq:laguerre_explicit}
 L_n^{(\alpha)}(t) = \sum_{j=0}^n 
\frac{\Gamma(n+\alpha+1)}{\Gamma(j+\alpha+1)(n-j)!}
\frac{(-t)^j}{j!}.
  \end{equation}
In terms of hypergeometric functions there is the following expression
\cite[formula 22.5.54]{abr:hmf}:
\begin{equation}
 \label{eq:laguerre_hyper}
  L_n^{(\alpha)}(t) = \frac{\Gamma(\alpha+n+1)}{\Gamma(\alpha+1)n!}
  \onefone{-n}{\alpha+1}t.
\end{equation}
Applying Kummer's transformation \cite[formula 13.1.27]{abr:hmf} 
to \eqref{eq:laguerre_hyper} leads to an alternative expression,
\begin{equation}
  \label{eq:laguerre_other}
  L_n^{(\alpha)}(t) = \rme^t \sum_{j=0}^\infty 
 \frac{\Gamma(\alpha+j+n)}{\Gamma(\alpha+j+1)n!}
  \frac{(-t)^j}{j!}.
\end{equation}

We denote by $\curlyW(g_1,\ldots,g_n)(x)$ the Wronskian
of the $n$ functions $g_1,\ldots,g_n$, evaluated at $x$:
\begin{equation}
  \label{eq:12}
  \curlyW(g_1,\ldots,g_n) \coloneq \det\left(
    \begin{array}{cccc}
      g_1 & g_2 & \cdots & g_n \\
      g_1' & g_2' & \cdots & g_n' \\
      \vdots & \vdots & \ddots & \vdots \\
      g_1^{(n-1)} & g_2^{(n-1)} & \cdots & g_n^{(n-1)} 
    \end{array}\right).
\end{equation}
We have,
\begin{proposition} \label{prop:laguerre_eval}
  For $k\in\N$ and $\zeta\in\R$,
  \begin{equation}
    \label{eq:13}
    H(k,\zeta) = (-1)^{k(k-1)/2} \frac{(2\pi)^N N!}{2^{2kN+N^2}}
      \rme^{-N|\zeta|} \curlyW(L_N^{(k)},L_{N+1}^{(k)},\ldots,
      L_{N+k-1}^{(k)})(-2|\zeta|).
  \end{equation}
\end{proposition}
\dimostrazione 
The integral on the right-hand side of \eqref{eq:5} is the averaged
moment of the characteristic polynomial of random matrices from the
Laguerre unitary ensemble. Such averages were considered by
\name{Br\'ezin} and \name{Hikami}, who \shew ed \cite[page 114]{bre:cpo} 
that
\begin{equation}
  \label{eq:43}
  \int_0^\infty\!\cdots\!\int_0^\infty \prod_{j=1}^N (t-y_j)^k y_j^k
\rme^{-y_j} \Delta(\vec{y})^2\,\rmd \vec{y} =
\frac{N!\prod_{j=0}^{N-1} c_j}{\prod_{\ell=0}^{k-1} \ell!}
\curlyW(p_N,\ldots,p_{N+k-1})(t),
\end{equation}
where $p_j=(-1)^j j! L_j^{(k)}$, and
\begin{align}
  \label{eq:44}
  c_j &= \int_0^\infty p_j(y)^2 y^k \rme^{-y}\,\rmd y \nonumber \\
  &=j!^2k!\binomial{j+k}k = j!(j+k)!
\end{align}
(the last line being a classical result of Laguerre polynomials
\cite{sze:op}).
We have
\begin{align}
 \nonumber
 \curlyW(p_N,\ldots,p_{N+k-1}) &= (-1)^{kN+k(k-1)/2
}\prod_{j=0}^{k-1}(N+j)!\,
  \curlyW(L_N^{(k)},\ldots,L_{N+k-1}^{(k)})\\
 &= (-1)^{kN+k(k-1)/2}\prod_{j=N-k}^{N-1}(k+j)!\,
  \curlyW(L_N^{(k)},\ldots,L_{N+k-1}^{(k)})
  \label{eq:45}
\end{align}
and
\begin{align}
  \nonumber
  \frac{\prod_{j=0}^{N-1} c_j}{\prod_{\ell=0}^{k-1}\ell!} &=
  \prod_{j=k}^{N-1} j! \prod_{j=0}^{N-1} (j+k)! \\
  \label{eq:46}
  &=\prod_{j=0}^{N-k-1} (j+k)!\prod_{j=0}^{N-1} (j+k)!
\end{align}
Putting \eqref{eq:45} and \eqref{eq:46} into \eqref{eq:43}, and setting
$t=-2\zeta$ gives
\begin{multline}
  \label{eq:47}
  \int_0^\infty\!\cdots\!\int_0^\infty \prod_{j=1}^N (y_j+2\zeta)^k y_j^k
\rme^{-y_j} \Delta(\vec{y})^2\,\rmd \vec{y} \\
= (-1)^{k(k-1)/2}N!\prod_{j=0}^{N-1} (j+k)!^2  
\curlyW(L_N^{(k)},L_{N+1}^{(k)},\ldots, L_{N+k-1}^{(k)})(-2\zeta).
\end{multline}
Together with proposition \ref{prop:integral}, this proves \eqref{eq:13}
for $\zeta>0$. For general $\zeta\in\R$, we note that the function 
$H(k,\zeta)$ defined by \eqref{eq:4} is an even continuous 
function of $\zeta$,
so that we may replace $\zeta$ by $|\zeta|$ wherever it occurs.
\finire

For $k=1$, the evaluation \eqref{eq:13} reduces to
\begin{equation}
  \label{eq:14}
  H(1,\zeta)=\frac{\pi^N N!}{2^{N^2+N}}\rme^{-N|\zeta|}
  L_N^{(1)}(-2|\zeta|),
\end{equation}
so that the integral \eqref{eq:4} is proportional to a single
Laguerre polynomial.  Although we were not able to find a use for 
\eqref{eq:13} for general $k$, we present in appendix \ref{app:b} an
elementary proof of \eqref{eq:conjecture}, based on \eqref{eq:14}, and
properties of Laguerre polynomials.

We are also able to write the Wronskian appearing in \eqref{eq:13} as
a Hankel determinant without derivatives, which may be of independent
interest.
\begin{proposition}
  If $L_N^{(k)}$ denotes a Laguerre polynomial, then
  \begin{equation}
    \label{eq:15}
    \curlyW(L_N^{(k)},\ldots,L_{N+k-1}^{(k)})(t) = \det\left(
      L_{N+k-1-(i+j)}^{(2k-1)}(t)\right)_{i,j=0,\ldots,k-1}.
  \end{equation}
\end{proposition}
\dimostrazione We make repeated use of the identity
\begin{equation}
  \label{eq:48}
  \frac{\rmd}{\rmd t} L_n^{(\alpha)}(t) = -L_{n-1}^{(\alpha+1)}(t),
\end{equation}
to get
\begin{align}
  \nonumber
   \curlyW(L_N^{(k)},\ldots,L_{N+k-1}^{(k)})(t) &= \det\left(
    \frac{\rmd^j}{\rmd t^j}  L_{N+i}^{(k)}(t)\right)_{i,j=0,\ldots,k-1} \\
    \nonumber
    &= \det\left( (-1)^j L_{N+i-j}^{(k+j)}(t)\right)_{i,j=0,\ldots,k-1} \\
   &=  (-1)^\omega\det\left(L_{N+i-j}^{(k+j)}(t)\right)_{i,j=0,\ldots,k-1}, 
\label{eq:49}
\end{align}
where $\omega=k/2$ if $k$ is even and $(k-1)/2$ if $k$ is odd. By
application of the identity 
\begin{equation}
  \label{eq:50}
  L_n^{(\alpha-1)}(t) + L_{n-1}^{(\alpha)}(t) = L_n^{(\alpha)}(t)
\end{equation}
{\it via\/} row operations on the matrix in \eqref{eq:49}, we get
\begin{equation}
  \label{eq:51}
    \curlyW(L_N^{(k)},\ldots,L_{N+k-1}^{(k)})(t)  = (-1)^\omega
    \det\left(L_{N+i-j}^{(2k-1)}(t)\right)_{i,j=0,\ldots,k-1}.
\end{equation}
Finally, switching rows according to $j\mapsto k-1-j$ gives $\omega$
transpositions, leading to \eqref{eq:15}.  \finire

\subsection{Evaluation in terms of hypergeometric functions of
matrix argument}\label{sec:eval-terms-hyperg}

The (single-variable) hypergeometric function $\vphantom{F}_p F_q$ is
defined formally \cite{bai:ghs} by the series
\begin{equation}
  \label{eq:30}
  \vphantom{F}_p F_q(a_1,\ldots,a_p; b_1,\ldots b_q; x) \coloneq 
  \sum_{j=0}^\infty
\frac{(a_1)_j\cdots(a_p)_j}{(b_1)_j\cdots (b_q)_j} \frac{x^j}{j!}
\end{equation}
where $(\cdot)_\cdot$ is the rising Pochhammer symbol. The
parameters $a_1,\ldots,a_p,b_1,\ldots,b_q$ can be arbitrary complex
numbers, however if $b_i \in \Z\smallsetminus\N$ then the series
\eqref{eq:30} becomes undefined unless there is a coresponding
parameter $a_{i'}\in \Z$ with $b_i<a_{i'}\leq 0$, in which case we
adopt the convention that the series \eqref{eq:30} terminates after
$a_{i'}$ terms.

The hypergeometric functions of a matrix argument provide a
multi-variable generalisation of \eqref{eq:30}. They have been studied
in \cite{kor:hti, kan:sia, yan:aco}, and have been found to occur in
the context of random matrix theory in the statistics of extreme
eigenvalues \cite{con:snd, kri:ote, dum:phd, dum:dot} and moments of
characteristic polynomials off the critical line \cite{for:dpe,
  for:sds}, amongst other places.
To generalise \eqref{eq:30}, the sum over integers in is replaced
by a sum over partitions, the Pochhammer symbols are replaced by
the generalised Pochhammer symbols defined in section 
\ref{sec:main-results}, and the univariate monomials $x^j$ are replaced
by Jack polynomials (see \cite{sta:scp} or \cite[chapter 12]{for:lga}).

Let $\sigma>0$ be a parameter, and $X$ be an $N\times N$ matrix
with eigenvalues $x_1,\ldots,x_N$. Then
\begin{align}
  \label{eq:31}
 \pFq (a_1,\ldots a_p; b_1,\ldots,b_q; X) &\coloneq
 \pFq (a_1,\ldots a_p; b_1,\ldots,b_q; x_1,\ldots,
x_N)\\
  \label{eq:32}
  &\coloneq \sum_{\lambda} \frac{[a_1]_\lambda^{(\sigma)}\cdots 
[a_p]_\lambda^{(\sigma)}}{[b_1]_\lambda^{(\sigma)}\cdots
[b_q]_\lambda^{(\sigma)}} \frac{C_\lambda^{(\sigma)}(x_1,\ldots,x_N)}
{|\lambda|!},
\end{align}
where the Jack polynomials $C_\lambda^{(\sigma)}$ are \cite{sta:scp}
homogeneous symmetric polynomial eigenfunctions of the partial
differential operator
\begin{equation}
  \label{eq:D_operator}
  \mathsf{D}(\sigma)\coloneq \frac{\sigma}2 \sum_{i=1}^N x_i^2
\frac{\partial^2}{\partial x_i^2} + \sum_{i\neq j} \frac{x_i^2}{x_i-x_j}
\frac{\partial}{\partial x_i},
\end{equation}
normalised so that
\begin{equation}
  \label{eq:33}
  C_\lambda^{(\sigma)}(1,\ldots,1) = \frac{\sigma^{2|\lambda|}|\lambda|!}
  {c_\lambda(\sigma)c'_\lambda(\sigma)}\left[ \frac{N}\sigma
    \right]_\lambda^{(\sigma)},
\end{equation}
(this is different to the normalisation used in \cite{sta:scp}) where
$c_\lambda(\sigma)$ and $c'_\lambda(\sigma)$ generalise the hook lengths
defined in section \ref{sec:main-results}; their definitions are
\begin{equation}
  \label{eq:34}
  c_\lambda(\sigma)\coloneq \prod_{\Box\in\lambda} \left( \ell(\Box) +1
    + \sigma g(\Box)\right),
\end{equation}
and
\begin{equation}
  \label{eq:35}
  c'_\lambda(\sigma)\coloneq \prod_{\Box\in\lambda} \left( \ell(\Box) +
\sigma(1   +  g(\Box))\right).
\end{equation}
For the case $\sigma=1$ the polynomials $C_\lambda^{(1)}$ are
proportional to Schur polynomials. In that case, the denominator
in \eqref{eq:33} is $c_\lambda(1)c'_\lambda(1) = h_\lambda^2$.

The convergence of the series in \eqref{eq:30} and \eqref{eq:32} depends
in general on the parameters $a_1,\ldots,a_p,b_1,\ldots,b_q$. If $p-q\leq1$
then the radii of convergence in \eqref{eq:30} and \eqref{eq:32} are
at least $1$.
If $a_i\in\Z\smallsetminus \N$ for any $i$ then the series terminate,
and are thus defined for all values of the arguments (in which case the
resulting functions are polynomials).

Hypergeometric functions of matrix argument enjoy a reflection
property \cite[page 812]{bor:Zmo} that we will make use of.
  Let $a\in\N$ and $b>a$. Then, 
  \begin{multline} \label{eq:hyper_reflect}
    \oneFone{\sigma}\! \left(-a;-b;x_1,\ldots,x_n\right) = \\
 ( x_1\cdots x_n )^a \prod_{j=1}^n \frac{\Gamma\! \left( b-a +
\frac{j-1}\sigma +1\right)}
{\Gamma\! \left( b+ \frac{j-1}\sigma +1 \right)}
\,\twoFzero{\sigma} \!\left( -a, 1+b-a+\frac{n-1}\sigma;;\frac{-1}{x_1},
\ldots, \frac{-1}{x_n}\right).
  \end{multline}
This may be viewed as a generalisation of the single-variable identity:
  \begin{equation} \label{eq:reflect}
  \onefone{-a}{-b}{x} = x^a\frac{\Gamma(b-a+1)}{\Gamma(b+1)}
    \vphantom{F}_2F_0\!\left( -a, 1+b-a ;; \frac{-1}x\right).
  \end{equation}
In both \eqref{eq:hyper_reflect} and \eqref{eq:reflect} we emphasise that
since the parameter $a$ is a negative integer, the hypergeometric
series are actually finite, and the left-hand side series
terminate before the denominators
in the summands in \eqref{eq:30} and \eqref{eq:32} become zero.
  
The main result of this section is the following:
  \begin{proposition} \label{prop:hyper_eval}
    For $k\in\N$ and $\zeta\in\R$,
    \begin{equation}
      \label{eq:16}
      H(k,\zeta) = \frac{\pi^N N!}{2^{(N+2k-1)N}}
      \tilde{F}_N(0,k) \rme^{-N|\zeta|} \,\oneFone{1}( -k; -2k; 2|\zeta|,
      \ldots,2|\zeta|).
    \end{equation}
  \end{proposition}

  The proof of proposition \ref{prop:hyper_eval} is based on the following
evaluation of the integral $L(k,\zeta)$:
\begin{proposition} \label{prop:hyper_eval2}
  Let $L$ denote the integral in \eqref{eq:6}, and let $k\in\N$ and
$\zeta>0$. Then
\begin{equation}
  \label{eq:17}
  L(k,\zeta) = \\
\Bigg(\prod_{j=0}^{N-1} \Gamma(j+2)\Gamma(k+1+j)\Bigg) (2\zeta)^{kN}
\twoFzero{1}\!\left( -k; N+k;; \frac{-1}{2\zeta},\ldots,\frac{-1}{2\zeta}
\right).
\end{equation}
\end{proposition}
\dimostrazione \name{Forrester} and \name{Keating} 
\cite[equation (3.2)]{for:sds} have 
proved\footnote{There is a typo in equation (3.2) of \cite{for:sds} in the
second parameter of the hypergeometric function. We have given here a
corrected formula.} the
following integral, valid for $\Re\{a\}>-1$, $\Re\{b\}>-1$ and either
$|t|>1$ or $2\mu\in\N$:
\begin{multline}
  \label{eq:36}
\frac1{S_N(a+1,b+1,\gamma)} \int_0^1\!\cdots\!\int_0^1 
  \prod_{j=1}^N x_j^a (1-x_j)^b(t-x_j)^{2\mu} |\Delta(\vec{x})|^{2\gamma}\,
  \rmd\vec{x} \\
 = t^{2\mu N} \twoFone{1/\gamma}\!\left(-2\mu, \gamma(N-1)+a+1;
  2\gamma(N-1)+a+b+2;\frac1t,\ldots,\frac1t\right),
\end{multline}
where
\begin{equation}
  \label{eq:37}
  S_N(\alpha,\beta,\gamma)\coloneq \prod_{j=0}^{N-1} \frac{\Gamma(\alpha+
j\gamma)\Gamma(\beta+j\gamma)\Gamma(1+(j+1)\gamma)}
{\Gamma(\alpha+\beta+(N+j-1)\gamma)\Gamma(1+\gamma)}.
\end{equation}

We let $2\mu=k\in\N$ and $a=k$, $b=L$, $\gamma=1$, $t=-2\zeta/L$, and 
make the changes of variables $x_j=y_j/L$. The integral on the left-hand
side of \eqref{eq:36} becomes 
\begin{equation}
  \label{eq:38}
  \frac{(-1)^{kN}}{L^{N^2+2kN}}\int_0^L\!\cdots\!\int_0^L
    \prod_{j=1}^N y_j^k\left(1-\frac{y_j}L\right)^L (y_j+2\zeta)^k
    \Delta(\vec{y})^2\,\rmd\vec{y}.
\end{equation}
Since $(1-y_j/L)^L\to\rme^{-y_j}$ as $L\to\infty$, and
\begin{equation}
  \label{eq:39}
  \left(1-\frac{y_j}L\right)^L\I_{[0,L]}(y) \leq \rme^{-y_j},
\end{equation}
for $y_j\geq 0$, we get, by the dominated convergence theorem, 
\begin{multline}
  \label{eq:40}
  L(k,\zeta)=\lim_{L\to\infty} (2\zeta)^{kN} \\ \times
L^{N^2+kN} S_N(k+1,L+1,1)
 \,\twoFone{1}\!
\left( -k, N+k; 2N+k+L; \frac{-L}{2\zeta},\ldots,
    \frac{-L}{2\zeta}\right),
\end{multline}
where the hypergeometric function of a matrix argument is a multivariate
polynomial since $k\in\N$.
By Stirling's formula and \eqref{eq:37} it follows that
\begin{equation}
  \label{eq:41}
  \lim_{L\to\infty} L^{N^2+kN}S_N(k+1,L+1,1) = \prod_{j=0}^{N-1}
  \Gamma(j+2)\Gamma(k+j+1).
\end{equation}
For a partition $\lambda$, it is clear that $[2N+k+L]_{\lambda}\sim
L^{|\lambda|}$ as $L\to\infty$, and the fact that $C_{\lambda}^{(1)}$
is homogeneous of degree $|\lambda|$ suffices to conclude from 
\eqref{eq:32} that
\begin{equation}
  \label{eq:42}
  \lim_{L\to\infty}\,   
\twoFone{1}\!\left( -k, N+k; 2N+k+L; \frac{-L}{2\zeta},\ldots,
    \frac{-L}{2\zeta}\right) = \twoFzero{1}\!\left(-k, N+k;; 
\frac{-1}{2\zeta},\ldots,\frac{-1}{2\zeta}\right).
\end{equation}
This completes the proof.
\finire

Similar integrals have been evaluated in \cite{bor:Zmo} by a different
method.

\dimostrazionea{proposition \ref{prop:hyper_eval}}
Combining the result of the proposition \ref{prop:hyper_eval2} with
the reflection formula \eqref{eq:hyper_reflect} and equation
\eqref{eq:5} leads to \eqref{eq:16}, for $\zeta>0$.
To pass to the case $\zeta\in\R$ we use the same argument as at the end
of the proof of proposition \ref{prop:laguerre_eval}. \finire

\mathversion{bold}
\section{Some integrals involving 
\texorpdfstring{$K_n(\epsilon,\zeta)$}{K_n(epsilon,zeta)}}
\mathversion{normal}
\label{sec:some-integr-involv}

In this section we will consider for integer values of $n$ and $p$,
integrals of the form
\begin{equation}
  \label{eq:23}
  \int_0^\infty K_n(\epsilon,\zeta)\rme^{-N\zeta}\zeta^p\,\rmd\zeta
\end{equation}
in the asymptotic r\'egime $\epsilon\downarrow 0$.

The function $K_n$ is defined in \eqref{eq:def_K} as a partial derivative
with respect to $\epsilon$ of a function of $\epsilon$ and $\zeta$. It
will be convenient to obtain alternative formul\ae\ for $K_n$ with
derivatives with respect to $\zeta$.
\begin{lemma} \label{lem:K_n} Let $K_n(\epsilon,\zeta)$ be defined by
  \eqref{eq:def_K}. If $n=2m$ is even then
\begin{equation}
  \label{eq:21}
  K_n(\epsilon,\zeta) = \frac{(-1)^m}{\pi} \frac{\partial^n}{\partial 
\zeta^n}\left( \frac{\epsilon}{\epsilon^2 + \zeta^2} \right),
\end{equation}
whereas if $n=2m-1$ is odd then
\begin{equation}
  \label{eq:21a}
  K_n(\epsilon,\zeta) = \frac{(-1)^{m+1}}{\pi} \frac{\partial^n}{\partial 
\zeta^n}\left( \frac{\zeta}{\epsilon^2 + \zeta^2} \right).
\end{equation}
\end{lemma}
\dimostrazione We begin by observing that
\begin{align}
  \label{eq:62}
  K_n(\epsilon,\zeta) &= \frac{(-1)^n}\pi \frac{\partial^n}{\partial
    \epsilon^n}\left( \frac\epsilon{\epsilon^2+\zeta^2}\right) \nonumber\\
  &= \frac{(-1)^n}{2\pi} \frac{\partial^n}{\partial \epsilon^n}
  \left( \frac1{\epsilon+\rmi\zeta} + \frac1{\epsilon-
\rmi\zeta}\right) \nonumber \\
  &= \frac{n!}{2\pi}\left( \frac1{\rmi^{n+1}(\zeta-\rmi\epsilon)^{n+1}} +
    \frac1{(-\rmi)^{n+1}(\zeta+\rmi\epsilon)^{n+1}}\right).
\end{align}
If $n=2m-1$ is odd, then
\begin{align}
  \label{eq:63}
  K_n(\epsilon,\zeta) &= \frac{n!}{2\pi \rmi^{2m}}\left(
    \frac1{(\zeta-\rmi\epsilon)^{n+1}} + \frac1{(\zeta+\rmi\epsilon)^{n+1}}
    \right) \nonumber \\
    &= \frac{(-1)^{n+m}}{2\pi} \frac{\partial^n}{\partial \zeta^n}
    \left( \frac1{\zeta-\rmi\epsilon} + \frac1{\zeta +\rmi\epsilon}\right)
    \nonumber \\
    &= \frac{(-1)^{m+1}}\pi \frac{\partial^n}{\partial\zeta^n}
    \left( \frac\zeta{\epsilon^2+\zeta^2}\right).
\end{align}
On the other hand, if $n=2m$ is even, then from \eqref{eq:62} we get
\begin{align}
  \label{eq:64}
  K_n(\epsilon,\zeta) &= \frac{(-1)^{n+m}}{2\pi\rmi} \frac{\partial^n}
  {\partial\zeta^n} \left( \frac1{\zeta-\rmi\epsilon} - \frac1
    {\zeta+\rmi\epsilon}\right) \nonumber \\
  &= \frac{(-1)^m}{\pi}\frac{\partial^n}{\partial \zeta^n}\left(
    \frac\epsilon{\epsilon^2+\zeta^2}\right).
\end{align}
\finire

\begin{lemma} \label{lem:diff_exp}
  Let $L_n^{(\alpha)}$ denote the Laguerre polynomials,
 defined by \eqref{eq:laguerre}. If $p\geq n$ then
\begin{equation}
   \label{eq:20}
   \frac{\partial^n}{\partial\zeta^n} (\zeta^p\rme^{-N\zeta}) =
   n! \zeta^{p-n}\rme^{-N\zeta}L_n^{(p-n)}(N\zeta).
 \end{equation}
For $p\leq n$, we have
\begin{equation}
  \label{eq:22}
   \frac{\partial^n}{\partial\zeta^n} (\zeta^p\rme^{-N\zeta}) 
= p! (-N)^{n-p}\rme^{-N\zeta}L_p^{(n-p)}(N\zeta).  
\end{equation}
\end{lemma}

\dimostrazione The case $p\geq n$ is slightly the simpler and we consider
it first. We have
\begin{align}
  \label{eq:65}
  \frac{\partial^n}{\partial \zeta^n} \left(\zeta^p\rme^{-N\zeta}
    \right) &= N^{-p} \frac{\partial^n}{\partial \zeta^n} \left(
      (N\zeta)^p \rme^{-N\zeta}\right) \nonumber \\
    &= N^{n-p} \left.\frac{\partial^n}{\partial\xi^n} \left(
        \xi^p \rme^{-\xi}\right)\right|_{\xi=N\zeta} \nonumber \\
    &= N^{n-p} n! \left. \xi^{p-n}\rme^{-\xi}L_n^{(p-n)}(\xi)
      \right|_{\xi=N\zeta} \nonumber \\
    &= n!\zeta^{p-n}\rme^{-N\zeta} L_n^{(p-n)}(N\zeta).
\end{align}
If $p\leq n$, we still have
\begin{equation}
  \label{eq:66}
  \frac{\partial^n}{\partial\zeta^n} (\zeta^p\rme^{-N\zeta}) =
  N^{n-p}\left.\frac{\partial^n}{\partial\zeta^n}(\xi^p\rme^{-\xi})
   \right|_{\xi=N\zeta},
\end{equation}
but now we write
\begin{align}
  \label{eq:67}
  \frac{\partial^n}{\partial\zeta^n} (\zeta^p\rme^{-N\zeta}) &=
  N^{n-p}\left.\frac{\partial^{n-p}}{\partial\xi^{n-p}}
    \left( \frac{\partial^p}{\partial\xi^p}(\xi^p\rme^{-\xi}
      )\right)\right|_{\xi=N\zeta} \nonumber \\
  &= N^{n-p} p!\left.\frac{\partial^{n-p}}{\partial \xi^{n-p}}
  \left( \rme^{-\xi}L_p^{(0)}(\xi)\right)\right|_{\xi=N\zeta}.
\end{align}
At this point we use \eqref{eq:laguerre_other} to write, for 
$\alpha\in\N_0$,
\begin{align}
  \label{eq:68}
  L_p^{(\alpha)}(\xi) = \rme^{\xi} \sum_{\ell=0}^\infty 
 \frac{(\alpha+p+\ell-1)!}{(\alpha+\ell)!p!}
 \frac{(-\xi)^\ell}{\ell!},
\end{align}
where the power series converges for all values of $\xi$. Thereby we get
an expression for $\rme^{-\xi}L_p^{(0)}(\xi)$ which
we may legitimately differentiate term-by-term to get
\begin{align}
  \label{eq:69}
  \frac{\partial^{n-p}}{\partial\xi^{n-p}}\left(
    \rme^{-\xi}L_p^{(0)}(\xi)\right) &= \sum_{\ell=n-p}^\infty
  \frac{(p+\ell-1)!}{\ell!p!}
 \frac{(-1)^\ell}{\ell!} \frac{\ell!}{(\ell-n+p)!}
  \xi^{\ell-n+p} \nonumber \\
  &= (-1)^{n-p}\sum_{m=0}^\infty \frac{(m+n-1)!}{(m+n-p)!p!}
  \frac{(-\xi)^m}{m!}, \qquad\mbox{{\it via\/} $\ell=m+n-p$,} \nonumber\\
  &= (-1)^{n-p} \rme^{-\xi} L_p^{(n-p)}(\xi),
\end{align}
using \eqref{eq:68} once more. Substitution of \eqref{eq:69} into
\eqref{eq:67} completes the proof.
\finire

For the case $p>n$ in \eqref{eq:23}, the asymptotic evaluations are
given by the following proposition.
\begin{proposition}\label{prop:K_drei}
  Let $p\geq n+1$. Then, as $\epsilon\downarrow 0$, we have
  \begin{equation}
    \label{eq:24}
    \int_0^\infty K_n(\epsilon,\zeta)\rme^{-N\zeta}\zeta^p \,\rmd\zeta
= \left\{\begin{array}{ll}
{\displaystyle \frac{n!(p-n-1)!}{\pi N^{p-n}}(-1)^{(n+1)/2} + \littleo(1),}
& \mbox{$n$ odd,}\\
\vphantom{\displaystyle\int} 
\littleo(1), & \mbox{$n$ even.}
\end{array}\right.
\end{equation}
\end{proposition}
\dimostrazione We use the representations of $K_n(\epsilon,\zeta)$
derived in lemma \ref{lem:K_n}, and then integrate by parts,
and insert the derivative formul\ae\ from lemma \ref{lem:diff_exp}.

If $n=2m-1$ is odd then
\begin{align}
 \int_0^\infty K_n(\epsilon,\zeta)\rme^{-N\zeta}\zeta^p\,\rmd\zeta
&= \frac{(-1)^{m+1}}\pi \int_0^\infty \frac{\partial^n}{\partial\zeta^n}
\left(\frac\zeta{\epsilon^2+\zeta^2}\right)\rme^{-N\zeta}\zeta^p\,
\rmd\zeta \nonumber \\
&= \frac{(-1)^{n+m+1}}\pi \int_0^\infty \frac\zeta{\epsilon^2+\zeta^2}
\frac{\partial^n}{\partial\zeta^n}(\rme^{-N\zeta}\zeta^p)\,\rmd\zeta
\nonumber \\
&= \frac{(-1)^m}\pi n!\int_0^\infty \frac{\zeta^{p-n+1}}{\epsilon^2
+\zeta^2}\rme^{-N\zeta}L_n^{(p-n)}(N\zeta)\,\rmd\zeta.   \label{eq:70}
\end{align}
Since $p-n+1\geq 2$ and since $|\zeta^2/(\epsilon^2+\zeta^2)|\leq 1$
the integral in \eqref{eq:70} is uniformly convergent in $\epsilon$,
so we can pass the limit under the integral and get
\begin{equation}
  \label{eq:72}
\lim_{\epsilon\downarrow0}
  \int_0^\infty K_n(\epsilon,\zeta)\rme^{-N\zeta}\zeta^p\,\rmd\zeta
= \frac{(-1)^m}\pi n! \int_0^\infty \zeta^{p-n-1}\rme^{-N\zeta}
L_n^{(p-n)}(N\zeta)\,\rmd\zeta.
\end{equation}
At this point, we insert the expansion \eqref{eq:laguerre_explicit}
for the Laguerre polynomial, to find
\begin{align}
\lim_{\epsilon\downarrow0}
    \int_0^\infty K_n(\epsilon,\zeta)\rme^{-N\zeta}\zeta^p\,\rmd\zeta
&= \frac{(-1)^m}\pi n! \sum_{\ell=0}^n \binomial{p}{n-\ell}
\frac{(-1)^\ell}{\ell!}N^\ell \int_0^\infty \zeta^{p-n+\ell-1}
\rme^{-N\zeta}\,\rmd\zeta \nonumber \\
&= \frac{(-1)^m}\pi \frac{n!(p-n-1)!}{N^{p-n}}
\sum_{\ell=0}^n (-1)^\ell \binomial{p}{n-\ell}\binomial{p-n+\ell-1}\ell
  \label{eq:71} \\
&= \frac{(-1)^m}\pi \frac{n!(p-n-1)!}{N^{p-n}},
\end{align}
using that the sum in \eqref{eq:71} evaluates to 1, a fact which is
proved in appendix \ref{app:b_1}.

In the case that $n=2m$ is even we again use lemma \ref{lem:K_n} and
lemma \ref{lem:diff_exp} to find that
\begin{align}
  \label{eq:78}
  \int_0^\infty K_n(\epsilon,\zeta)\rme^{-N\zeta}\zeta^p\,\rmd\zeta &=
  \frac{(-1)^{m+n}}\pi \int_0^\infty \frac{\epsilon}{\epsilon^2+\zeta^2}
  \frac{\partial^n}{\partial\zeta^n}(\rme^{-N\zeta}\zeta^p)\,\rmd\zeta
  \nonumber \\
  &=\frac{(-1)^{m+n}}\pi n!\int_0^\infty \frac{\epsilon}{\epsilon^2
    +\zeta^2} \zeta^{p-n}\rme^{-N\zeta} L_n^{(p-n)}(N\zeta)\,\rmd\zeta
  \nonumber \\
  &=\littleo(1),
\end{align}
as $\epsilon\downarrow0$, using lemma \ref{lem:delta} below.
\finire

In the proof of propositon \ref{prop:K_drei} we used the following
standard result (quoted without proof):
\begin{lemma} \label{lem:delta}
  If $f$ is a bounded function with $f(x)\to f_0$ as $x\downarrow 0$,
then
\begin{equation}
  \label{eq:28}
  \int_0^{\infty}\frac{\epsilon}{\epsilon^2 + x^2}f(x)\,\rmd x =
   \frac{\pi f_0}2 + \littleo(1),
\end{equation}
as $\epsilon\downarrow 0$.
\end{lemma}

\begin{proposition} \label{prop:K_vier}
  If $p\leq n$ and $n=2m$ is even, then
  \begin{equation}
    \label{eq:25}
    \int_0^\infty K_n(\epsilon,\zeta)\rme^{-N\zeta} \zeta^p\,\rmd \zeta
    = \frac{(-1)^{n/2}}2 \frac{n!}{(n-p)!}(-N)^{n-p} + \littleo(1),
  \end{equation}
as $\epsilon\downarrow 0$.
\end{proposition}
\dimostrazione Applying lemma \ref{lem:K_n} and lemma
\ref{lem:diff_exp} we get
\begin{align}
  \label{eq:79}
  \int_0^\infty K_n(\epsilon,\zeta)\rme^{-N\zeta}\zeta^p\,\rmd\zeta &=
  \frac{(-1)^{m+n}}\pi \int_0^\infty \frac\epsilon{\epsilon^2+\zeta^2}
  \frac{\partial^n}{\partial\zeta^n}(\rme^{-N\zeta}\zeta^p)\,\rmd\zeta
  \nonumber \\
  &=\frac{(-1)^{m+n}}\pi p!(-N)^{n-p}\int_0^\infty
  \frac\epsilon{\epsilon^2+\zeta^2}\rme^{-N\zeta}L_p^{(n-p)}(N\zeta)\,
  \rmd\zeta.
\end{align}
Applying lemma \ref{lem:delta}, we find in the limit
$\epsilon\downarrow0$,
\begin{align}
  \label{eq:80}
  \int_0^\infty K_n(\epsilon,\zeta)\rme^{-N\zeta}\zeta^p\,\rmd\zeta &=
  \frac{(-1)^{m+n}}\pi p!(-N)^{n-p}\frac\pi2 L_n^{(n-p)}(0) +
  \littleo(1) \nonumber \\
  &=\frac{(-1)^{m+n}}2 \frac{n!}{(n-p)!}(-N)^{n-p} + \littleo(1),
\end{align}
using the explicit representation \eqref{eq:laguerre_explicit} for the
Laguerre polynomial. \finire

For odd $n\geq p$, we have the slightly more subtle result:
\begin{proposition} \label{prop:K_fuenf}
Let $f(\zeta)\coloneq \sum_{p=0}^P f_p \zeta^p$, with $P\leq n$, where
$n$ is now odd. Then provided that
\begin{equation}
  \label{eq:26}
  \sum_{p=0}^P \binomial{n}p \frac{f_p p!}{(-N)^p} = 0,
\end{equation}
then
\begin{equation}
  \label{eq:27}
  \lim_{\epsilon\downarrow 0}\int_0^{\infty} K_n(\epsilon,\zeta)
  \rme^{-N\zeta}f(\zeta)\,\rmd\zeta = \frac{(-1)^{(n+1)/2}}\pi
    \sum_{p=1}^P\sum_{\ell=1}^p p!\binomial{n}{p-\ell}
    \frac{(-1)^\ell (-N)^{n-p}}{\ell} f_p.
\end{equation}
If condition \eqref{eq:26} does not hold, then the limit
in \eqref{eq:27} diverges.
\end{proposition}
\dimostrazione   Following the arguments above, we find for $n=2m-1$ and
$p\leq n$,
\begin{equation}
  \label{eq:81}
  \int_0^\infty K_n(\epsilon,\zeta)\rme^{-N\zeta}\zeta^p\,\rmd\zeta
= \frac{(-1)^{m}}\pi p!(-N)^{n-p}\int_0^\infty \frac\zeta
{\epsilon^2+\zeta^2}\rme^{-N\zeta}L_p^{(n-p)}(N\zeta)\,\rmd\zeta.
\end{equation}
Since 
\begin{equation}
  \label{eq:83}
  L_p^{(n-p)}(N\zeta) = \binomial{n}p + \sum_{\ell=1}^p
  \binomial{n}{p-\ell}\frac{(-N)^\ell}{\ell!}\zeta^\ell,
\end{equation}
we have that
\begin{align}
  \label{eq:84}
  \int_0^\infty \frac{\zeta\rme^{-N\zeta}}{\epsilon^2+\zeta^2}
  L_p^{(n-p)}(N\zeta)\,\rmd\zeta &= \binomial{n}p\int_0^\infty
  \frac{\zeta\rme^{-N\zeta}}{\epsilon^2+\zeta^2}\,\rmd\zeta
   \\
&\qquad + \sum_{\ell=1}^p \binomial{n}{p-\ell}\frac{(-N)^\ell}{\ell!}
   \int_0^\infty \zeta^{\ell-1}\rme^{-N\zeta}\,\rmd\zeta + \littleo(1),
\nonumber
\end{align}
passing the limit $\epsilon\downarrow0$ in a similar way as in the proof
of the first part of proposition \ref{prop:K_drei}.
Evaluating the integral on the right hand side of \eqref{eq:84} we
get
\begin{multline}
  \label{eq:110}
  \int_0^\infty K_n(\epsilon,\zeta)\rme^{-N\zeta}\zeta^p\,\rmd\zeta
= \frac{(-1)^m}\pi p!(-N)^{n-p}\binomial{n}{p}\int_0^\infty
\frac{\zeta\rme^{-N\zeta}}{\epsilon^2+\zeta^2}\,\rmd\zeta \\+
\frac{(-1)^m}\pi p!(-N)^{n-p}\sum_{\ell=1}^p\binomial{n}{p-\ell}
\frac{(-1)^\ell}\ell+\littleo(1).
\end{multline}
It is clear that when condition \eqref{eq:26} is satisfied, the 
contributions coming from the first term of the right-hand side of
\eqref{eq:110} cancel, and we arrive to \eqref{eq:27}.

If condition \eqref{eq:26} does not hold, then we need to prove that
the integral on the right-hand side of \eqref{eq:110} diverges
as $\epsilon\downarrow0$. To see this, we integrate by parts to get
\begin{equation}
  \label{eq:86}
  \int_0^\infty \frac{\zeta}{\epsilon^2+\zeta^2}\rme^{-N\zeta}\,\rmd\zeta
=-\log\epsilon+\frac{N}{2}\int_0^\infty\log(\epsilon^2+\zeta^2)
\rme^{-N\zeta}\,\rmd\zeta
\end{equation}
and observe that the integral on the right-hand side of \eqref{eq:86}
is $\Ord(1)$ as $\epsilon\downarrow0$. \finire

\section{Moments of characteristic polynomials and their derivatives}
\label{sec:sechs}
We have now collected ingredients required to prove theorem
\ref{thm:main}. After giving the proof below, we then will consider
the asymptotics limit $N\to\infty$ of large matrix size 
(proposition \ref{prop:main_limit} below).

\mathversion{bold}
\subsection{The finite $N$ case}
\mathversion{normal}

\dimostrazionea{theorem \ref{thm:main}} We use proposition \ref{prop:zwei}
and proposition \ref{prop:hyper_eval} to write, for $k\in\N$,
\begin{equation}
 \tilde{F}_N(h,k) = \lim_{\epsilon\downarrow 0}  \frac{1}{2^{2h}}
 \tilde{F}_N(0,k) \int_{-\infty}^\infty K_{2h}(\epsilon,\zeta)
\rme^{-N|\zeta|}\,\oneFone{1}(-k;-2k; 2|\zeta|,\ldots,2|\zeta|)\,\rmd\zeta.
  \label{eq:87}
\end{equation}
We shall derive an expansion for the hypergeometric function
in terms of the coefficients $\curlyC_N(p,k)$ defined by
\eqref{eq:curlyCN}. We begin by using the definition \eqref{eq:32} to write
\begin{equation}
  \label{eq:90}
   \oneFone{1}(-k;-2k;2|\zeta|,\ldots,2|\zeta|) =
 \sum_\lambda \frac{[-k]_\lambda}{[-2k]_\lambda} \frac{C_\lambda^{(1)}(
2|\zeta|,\ldots,2|\zeta|)}{|\lambda|!}.
\end{equation}
The series in \eqref{eq:90} is finite.  Indeed, the sum runs over only
those partitions $\lambda$ with largest part not greater than $k$ (so
that the factor $[-2k]_\lambda$ in the denominator is never zero) and
the Jack polynomial $C_\lambda^{(1)}$ vanishes if $\lambda$ is a
partition with more than $N$ parts. This places an upper bound of $kN$
on the sum of the parts of $\lambda$.  We use the homogeneity of Jack
polynomials, and the normalisation \eqref{eq:33} to get
\begin{align}
 \oneFone{1}(-k;-2k;2|\zeta|,\ldots,2|\zeta|) &=
 \sum_\lambda \frac{[-k]_\lambda}{[-2k]_\lambda} (2|\zeta|)^{|\lambda|}
\frac{C_\lambda^{(1)}(1,\ldots,1)}{|\lambda|!} \nonumber \\
 &= \sum_\lambda \frac{[-k]_\lambda}{[-2k]_\lambda} (2|\zeta|)^{|\lambda|}
\frac{[N]_\lambda}{h_\lambda^2}.
 \label{eq:88}
\end{align}
We can index the sum in \eqref{eq:88} by transposes of partitions,
rather than the partitions themselves. This gives
\begin{align}
  \oneFone{1}(-k;-2k;2|\zeta|,\ldots,2|\zeta|) &=
\sum_\lambda \frac{[-k]_{\lambda^{\rm T}}}{[-2k]_{\lambda^{\rm T}}} 
\frac{[N]_{\lambda^{\rm T}}}{h_{\lambda^{\rm T}}^2}
(2|\zeta|)^{|\lambda|} \nonumber \\
&= \sum_\lambda \frac{[k]_\lambda (-1)^{|\lambda|}[-N]_\lambda}
{[2k]_\lambda h_\lambda^2} (2|\zeta|)^{|\lambda|},
  \label{eq:107}
\end{align}
using \eqref{eq:106}. We will group the terms of \eqref{eq:107} so
that partitions of the same integer $p$ are summed together.  Since we
transposed the partitions, we know that $\lambda$ can have at most $k$
parts.  This manipulation brings us finally to
\begin{align}
    \oneFone{1}(-k;-2k;2|\zeta|,\ldots,2|\zeta|) &=
\sum_{p=0}^{kN} \Bigg( \sum_{\lambda\vdash_k p} 
\frac{[k]_\lambda [-N]_\lambda}
{[2k]_\lambda h_\lambda^2}\Bigg) (-2|\zeta|)^p \nonumber \\
 &= \sum_{p=0}^{kN} \curlyC_N(p,k)|\zeta|^p.
 \label{eq:108}
\end{align}

Substituting \eqref{eq:108} into \eqref{eq:87} we arrive at
\begin{equation}
  \label{eq:89}
 \tilde{F}_N(h,k) = \lim_{\epsilon\downarrow 0}  \frac{2}{2^{2h}}
 \tilde{F}_N(0,k) \sum_{p=0}^{kN} \curlyC_N(p,k)
 \int_{0}^\infty K_{2h}(\epsilon,\zeta)
 \rme^{-N\zeta}\zeta^p\,\rmd\zeta.
\end{equation}
We split the sum into two contributions according to $0\leq p \leq 2h$
and $2h < p \leq kN$, and apply proposition \ref{prop:K_drei} to the
second sum, and apply proposition \ref{prop:K_fuenf} to the first
sum. By proposition \ref{prop:zwei} we know {\it a priori\/} that the
limit $\epsilon\downarrow 0$ exists in \eqref{eq:89}, so condition
\eqref{eq:26} must hold with $f_p = \curlyC_N(p,k)$. This means that
in addition to proving \eqref{eq:3}, we have also proved
the combinatorial identity
\begin{equation}
  \label{eq:91}
  \sum_{p=0}^{2h} \binomial{2h}p \curlyC_N(p,k)\frac{p!}{(-N)^p} =0,
\end{equation}
valid for $2h$ an odd integer with $0< h\leq k$. \finire

We remark that starting from equation \eqref{eq:89} (which does not
depend on the parity of $2h$), and using propositions 
\ref{prop:K_drei} and \ref{prop:K_vier}, we can re-prove \name{Dehaye}'s
result, theorem \ref{thm:dehaye}, for $2h$ even and $k$ integer,
using our methods.

\mathversion{bold}
\subsection{The \texorpdfstring{$N\to\infty$}{N to infinity} limit}
\mathversion{normal}

In order to pass to the limit $N\to\infty$ in \eqref{eq:3} we require an
estimate on the size of the coefficients $\curlyC_N(p,k)$. This is 
provided by the following lemma.
\begin{lemma}
  \label{lem:curlyC}
  Let $N\geq1$ and $p\geq 2$. Then
  \begin{equation}
    \label{eq:92}
    \curlyC_N(p,k) = \Ord_k\!\left( \frac{N^p}{p!} \right),
  \end{equation}
where the implied constant may depend on $k$, but is independent
of $N$ and $p$.
\end{lemma}
\dimostrazione
We use the fact that if $\lambda=(\lambda_1,\ldots,\lambda_k)$ is a
partition of $p$ into not more than $k$ parts,
then $\lambda_1\geq\lfloor\frac{p}k\rfloor$, where $\lfloor\cdot\rfloor$
is the integer-part function. Then we have,
\begin{align}
[2k]_\lambda = \prod_{i=1}^k \prod_{j=1}^{\lambda_i} (2k+j-i)
 &\geq \prod_{j=1}^{\lambda_1} (2k+j-1) \nonumber \\
 &\geq \prod_{j=1}^{\lfloor p/k\rfloor} (2k+j-1) \nonumber \\
 &=\frac{\Gamma(2k+\lfloor p/k\rfloor)}{\Gamma(2k)}.  
  \label{eq:93}
\end{align}
With this inequality, and the trivial estimate
\begin{equation}
  \label{eq:95}
  |[-N]_\lambda| \leq (N+k)^p,
\end{equation}
we can bound $\curlyC_N(p,k)$ as follows:
\begin{align}
|\curlyC_N(p,k)| = 2^p\Bigg|\sum_{\lambda \vdash_k p} \frac{[k]_\lambda 
  [-N]_\lambda}{[2k]_\lambda h_\lambda^2}\Bigg| &\leq
  \frac{2^p(N+k)^p \Gamma(2k)}{\Gamma(2k+\lfloor p/k\rfloor)}
  \sum_{\lambda \vdash_k p} \frac{[k]_\lambda}{h_\lambda^2} \nonumber \\
  &=\frac{2^p(N+k)^p \Gamma(2k)}{\Gamma(2k+\lfloor p/k\rfloor)}
  \sum_{\lambda \vdash p} \frac{[k]_\lambda}{h_\lambda^2},
 \label{eq:82}
\end{align}
where the last equality holds, since $[k]_\lambda=0$ if $\lambda$ is
a partition with more than $k$ parts. It follows from the
hook-content formula (see Theorem 7.21.2 of \cite{sta:ec2}
combined with Proposition 2.2 of \cite{han:sca})
that
\begin{equation}
  \label{eq:96}
  \sum_{\lambda\vdash p} \frac{[k]_\lambda}{h_\lambda^2} = \frac{k^p}
  {p!},   
\end{equation}
so we have proved 
\begin{equation}
  \label{eq:97}
  |\curlyC_N(p,k)|\leq \frac{N^p}{p!}\left(1+\frac{k}N\right)^p
  \frac{(2k)^p\Gamma(2k)}{\Gamma(2k+\lfloor p/k \rfloor)},
\end{equation}
which furnishes the required estimate. \finire

Lemma \ref{lem:curlyC} is probably far from optimal, but is sufficient to 
prove the result following.
\begin{proposition}
  \label{prop:main_limit}
  Let $h=(2m-1)/2$ for $m\in\N$ and let $k\in\N$ with $k>h-\frac12$.
  Then
  \begin{multline}
    \label{eq:NtoInf}
    \tilde{F}(h,k) = \frac{2(-1)^{h+1/2}}{2^{2h}\pi} 
  \tilde{F}(0,k)\Bigg\{
     \sum_{p=1}^{2h}\sum_{\ell=1}^p
 \binomial{2h}{p-\ell} \frac{(-1)^{\ell+2h-p}}\ell
 p!\, \curlyC(p,k) + \\ \sum_{p=2h+1}^\infty
 (2h)!(p-2h-1)!\, \curlyC(p,k)\Bigg\},
  \end{multline}
where
\begin{equation}
  \label{eq:102}
  \tilde{F}(0,k) = \prod_{j=1}^k \frac{\Gamma(j)}{\Gamma(k+j)}.
\end{equation}
\end{proposition}
\dimostrazione
We recall that 
\begin{equation}
  \label{eq:98}
  \tilde{F}(h,k) \coloneq \lim_{N\to\infty} \frac1{N^{k^2+2h}}
\tilde{F}_N(h,k).
\end{equation}
We can apply the limit term-by-term in the first sum of \eqref{eq:3},
using \eqref{eq:curlyC_asympt} and the fact that
\begin{equation}
  \label{eq:103}
  \tilde{F}_N(0,k) \sim \tilde{F}(0,k) N^{k^2},
\end{equation}
as $N\to\infty$, which was proved in \cite{kea:rmtz}.

For the second sum, lemma \ref{lem:curlyC} gives us 
\begin{equation}
  \label{eq:100}
  \frac{(p-2h-1)!}{N^p}\curlyC_N(p,k) = \Ord_k\!\left(\frac1{p^{2h+1}}
    \right),
\end{equation}
so that the second summand of \eqref{eq:3} is bounded independently
of $N$ by a summable function of $p$. Taken together with 
\eqref{eq:curlyC_asympt}, this allows us to
apply Tannery's theorem \cite[\S 49]{bro:ait} to prove
\begin{equation}
  \label{eq:111}
  \lim_{N\to\infty}\sum_{p=2h+1}^{kN} \frac{(p-2h-1)!}{N^p}
\curlyC_N(p,k) = \sum_{p=2h+1}^\infty (p-2h-1)!\,\curlyC(p,k),
\end{equation}
and hence we get \eqref{eq:NtoInf}. \finire
\section{Partition sums and a proof of (\ref{eq:conjecture})}
\label{sec:partition-sums-proof}

In order to give explicit formul\ae\ for the moment $\tilde{F}(h,k)$
we require closed forms for $\curlyC(p,k)$. We have been able to find
these forms for $k=1$ and $k=2$:
\begin{proposition} \label{prop:curlyC_evals}
  Let $p\in\N$ and $\curlyC(p,k)$ be defined by \eqref{eq:curlyC}. Then,
  \begin{equation}
    \label{eq:18}
    \curlyC(p,1) = \frac{2^p}{p!(p+1)!},
  \end{equation}
and
\begin{equation}
  \label{eq:19}
  \curlyC(p,2) = \frac{12(2p+4)!2^p}{p!(p+2)!(p+3)!(p+4)!}.
\end{equation}
\end{proposition}
\dimostrazione We recall that
\begin{equation*}
\curlyC(p,k) \coloneq 2^p \sum_{\lambda\vdash_k p} \frac{[k]_{\lambda}}
{[2k]_\lambda h_{\lambda}^2},
\end{equation*}
where the summation goes over partitions of $p$ into not more than
$k$ parts. In the case $k=1$, then only the single partiton $\lambda=(p)$
is admitted. In this case, it is easy to see that $[1]_\lambda=p!$,
$[2]_\lambda=(p+1)!$ and $h_\lambda=p!$. This immediately leads to
\eqref{eq:18}.

\begin{figure}
  \centering
  \begin{tikzpicture}
   \draw (0.5,0) -- (0.5,1);
    \draw (1,0) -- (1,1);
    \draw (2,0) -- (2,1);
    \draw (3.2,0) -- (3.2,1);
    \draw (3.7,0.5) -- (3.7,1);
    \draw (4.9,0.5) -- (4.9,1);
    \draw (5.9,0.5) -- (5.9,1);
    \draw (-0.3,0) rectangle (3.2,0.5); 
    \draw (-0.3,0.5) rectangle (7.7,1);
   \path (0.1,0.25)  node  {$-1$}
          (0.1,0.75)  node  {$0$}
          (0.75,0.25) node  {$0$}
          (0.75,0.75) node  {$1$}
          (1.5,0.25) node  {$\cdots$}
          (1.5,0.75) node  {$\cdots$}
          (2.6,0.25) node  {$x-2$}
          (2.6,0.75) node  {$x-1$}
          (3.45,0.75) node {$x$}
          (4.3,0.75) node {$x+1$}
          (5.4,0.75) node {$\cdots$}
          (6.8,0.75) node {$p-x-1$};
   \draw[decorate,decoration=brace] (3.2,-0.2) -- (-0.3,-0.2) node [midway,
anchor=north] {$x$ boxes};
   \draw[decorate,decoration=brace] (-0.3,1.2) -- (7.7,1.2) node [midway,
anchor=south] {$p-x$ boxes};
  \end{tikzpicture}
  \caption{Partitions of $p$ into $2$ parts. The contents of each box are 
the values of $j-i$ appearing in the definition of the Pochhammer symbol.}
  \label{fig:partitions_evals_eins}
\end{figure}
\begin{figure}
  \centering
  \begin{tikzpicture}
    \draw (0,0) -- (0,1);
    \draw (1.8,0) -- (1.8,1);
    \draw (3,0) -- (3,1);
    \draw (4,0) -- (4,1);
    \draw (6,0) -- (6,1);
    \draw (7.4,0.5) -- (7.4,1);
    \draw (8.4,0.5) -- (8.4,1);
    \draw (8.9,0.5) -- (8.9,1);
    \draw (9.4,0.5) -- (9.4,1);
    \draw (0,0) -- (6,0);
    \draw (0,0.5) -- (9.4,0.5);
    \draw (0,1) -- (9.4,1);
    \path (0.9,0.75) node {$p-x+1$}
          (0.9,0.25) node {$x$}
          (2.4,0.75) node {$p-x$}
          (2.4,0.25) node {$x-1$}
          (3.5,0.75) node {$\cdots$}
          (3.5,0.25) node {$\cdots$}
          (5,0.75) node {$p-2x+2$}
          (5,0.25) node {$1$}
          (6.7,0.75) node {$p-2x$}
          (7.9,0.75) node {$\cdots$}
          (8.65,0.75) node {$2$}
          (9.15,0.75) node {$1$};
  \end{tikzpicture}
  \caption{The hooks for the partition $\lambda=(p-x,x)$. The hook
length $h_\lambda$ is the product of the entries in the boxes.}
  \label{fig:partitions_evals_zwei}
\end{figure}

For $k=2$ we require all partitions of $p$ into not more than $2$ parts.
These partitions are of the form $\lambda=(p-x,x)$, where $0\leq x <
(p+1)/2$ (see figure \ref{fig:partitions_evals_eins}). For a partition
of this form, we have
\begin{equation}
  \label{eq:52}
  [2]_\lambda = (p-x+1)!x!
\end{equation}
and
\begin{equation}
  \label{eq:53}
  [4]_\lambda = \frac{(x+2)!(p-x+3)!}{12},
\end{equation}
and the hook length $h_\lambda$ is given by (see 
figure \ref{fig:partitions_evals_zwei})
\begin{equation}
  \label{eq:54}
  h_\lambda = \frac{x!(p-x+1)!}{p-2x+1},
\end{equation}
so that
\begin{equation}
  \label{eq:55}
  \curlyC(p,2) = 2^p\sum_{0\leq x < (p+1)/2} \frac{12(p-2x+1)^2}{x!(x+2)!
    (p-x+3)!(p-x+1)!}.
\end{equation}
We observe that the summand in \eqref{eq:55} is invariant under the
reflection $x\mapsto p+1-x$, so that the half-range sum can
be replaced by one half times the sum from $0$ to $p+1$, giving
\begin{equation}
  \label{eq:56}
  \curlyC(p,2) = 2^p 6 \sum_{x=0}^{p+1} \frac{(p-2x+1)^2}{x!(x+2)!
    (p-x+3)!(p-x+1)!}.
\end{equation}
In appendix \ref{app:b_2} the sum in \eqref{eq:56} is evaluated, 
whereupon \eqref{eq:19} follows. 
\finire

Based on the results of proposition \ref{prop:curlyC_evals}, it is
tempting to conjecture that $\curlyC(p,k)$ will be $2^p$ times a ratio
of products of factorials for all $k\in\N$. However a computer-based
investigation of $\curlyC(p,3)$ has \shew n that this structure appears
to break down when $k=3$.

When $h=1/2$ and $k=1$, substituting \eqref{eq:18} into
equation \eqref{eq:NtoInf} gives
\begin{align}
  \tilde{F}({\textstyle\frac12},1) &= \frac1\pi\left( 1 - 
 \sum_{p=2}^\infty \frac{2^p (p-2)!}{p!(p+1)!}\right) \nonumber \\
&=\frac{\rme^2-5}{4\pi}, 
  \label{eq:94}
\end{align}
as conjectured.

Using \eqref{eq:19} for the case $k=2$ in \eqref{eq:NtoInf}, we can
evaluate\footnote{We have used the algebraic manipulation package {\tt
    Maple} to derive closed forms for the sums in terms of
  hypergeometric functions.}
\begin{align}
  \tilde{F}({\textstyle\frac12},2) &= \frac1{180\pi}\left(
 15 -  7\vphantom{F}_3F_3\left( 1,1,{\textstyle\frac92}; 3,6,7; 8\right)
    \right)\nonumber \\
  &\approx 0.008\,15\ldots
  \label{eq:104}
\end{align}
and
\begin{align}
  \tilde{F}({\textstyle\frac32},2) &= \frac1{10080\pi}\left(
    33\vphantom{F}_3F_3\left( 1,1,{\textstyle\frac{13}2}; 5,8,9; 8\right)
    -28 \right)\nonumber \\
  &\approx 0.000\,354\ldots
  \label{eq:105}
\end{align}
These last two values do not appear to have been derived or
conjectured before.

\subsection*{Acknowledgements}
I am grateful for a number of interesting conversations regarding
this work with 
\href{http://www.math.tamu.edu/~berko/}{Gregory \name{Berkolaiko}},
\href{http://maths.york.ac.uk/www/ch540}{Chris \name{Hughes}}, 
\href{http://www.math.ethz.ch/~pdehaye/}{Paul-Olivier \name{Dehaye}}, 
\href{http://www.lboro.ac.uk/departments/ma/people/hallnas.html}
{Martin \name{Halln\"as}} and
\href{http://www-staff.lboro.ac.uk/~mait/personal2/index.html}
{Ian \name{Thompson}}.
\href{http://www.maths.bris.ac.uk/people/profile/majpk}{Jon 
\name{Keating}} 
provided a number of crucial suggestions.
I am grateful to \href{http://www.ms.unimelb.edu.au/~matpjf/matpjf.html}
{Peter \name{Forrester}} for bringing reference \cite{bor:Zmo} to my 
attention.

\appendix

\mathversion{bold}
\section{An elementary evaluation of 
\texorpdfstring{$\tilde{F}(\frac12,1)$}{Ftwiddle(1/2,1)}}
\mathversion{normal}  \label{app:b}
In this appendix we give a completely elementary derivation of 
\eqref{eq:conjecture} based on classical properties of Laguerre
polynomials (see, for example, \cite{sze:op}).

We first recall the alternative expression \eqref{eq:21a} for 
$K_1(\epsilon,\zeta)$:
 \begin{equation} \label{eq:app_B:uno}
    K_1(\epsilon,\zeta) = \frac1{\pi} \frac{\partial}{\partial \zeta}
    \left( \frac{\zeta}{\epsilon^2+\zeta^2} \right).
  \end{equation}
Equipped with this, we proceed from the integral representation of
proposition \ref{prop:zwei} and equation \eqref{eq:14},
\begin{displaymath}
  \tilde{F}_N({\textstyle\frac12},1) = \lim_{\epsilon\downarrow 0} \frac12
\int_{-\infty}^\infty K_1(\epsilon,\zeta)\rme^{-N|\zeta|}
 L_N^{(1)}(-2|\zeta|)\,\rmd\zeta.
\end{displaymath}
Substituting \eqref{eq:app_B:uno} and integrating-by-parts, we get
\begin{align}
  \frac1\pi \int_0^\infty K_1(\epsilon,\zeta)\rme^{-N\zeta}
  L_N^{(1)}(-2\zeta)\,\rmd \zeta &=
  -\frac1\pi\int_0^\infty \frac{\zeta}{\epsilon^2+\zeta^2}
  \frac{\partial}{\partial\zeta}\left( \rme^{-N\zeta}
    L_N^{(1)}(-2\zeta)\right)\,\rmd \zeta \nonumber \\
  &=\frac1\pi\int_0^\infty \frac{\zeta}{\epsilon^2+\zeta^2}
  \rme^{-N\zeta}\left( NL_N^{(1)}(-2\zeta)-2L_{N-1}^{(2)}(-2\zeta)
    \right)\,\rmd\zeta, \nonumber
\end{align}
where we have used the fact that 
\begin{displaymath}
  \frac{\rmd}{\rmd t} L_N^{(\alpha)}(t) = -L_{N-1}^{(\alpha+1)}(t).
\end{displaymath}
By a standard recurrence for Laguerre polynomials,
\begin{displaymath}
  NL_N^{(1)}(-2\zeta) - 2L_{N-1}^{(2)}(-2\zeta) = 2\zeta 
L_{N-1}^{(3)}(-2\zeta),
\end{displaymath}
so that we have
\begin{displaymath}
  \tilde{F}_N({\textstyle\frac12},1) = \frac2\pi \lim_{\epsilon
    \downarrow 0} \int_0^{\infty} \frac{\zeta^2}{\epsilon^2+\zeta^2}
  \rme^{-N\zeta} L_{N-1}^{(3)}(-2\zeta)\,\rmd\zeta.
\end{displaymath}
Since
\begin{displaymath}
  0\leq \frac{\zeta^2}{\epsilon^2+\zeta^2}\leq 1,
\end{displaymath}
we can apply the dominated convergence theorem to pass the limit under
the integral, to get 
\begin{equation}
  \tilde{F}_N({\textstyle\frac12},1) = \frac2\pi \int_0^{\infty}
  \rme^{-N\zeta} L_{N-1}^{(3)}(-2\zeta)\,\rmd\zeta. 
\end{equation}
At this point, we insert \eqref{eq:laguerre_explicit}, the explicit 
representation of
$L_N^{(\alpha)}$, getting that
\begin{align*}
  \int_0^{\infty} \rme^{-N\zeta}L_{N-1}^{(3)}(-2\zeta)\,\rmd \zeta
  &= \sum_{n=0}^{N-1} \left(
    \begin{array}{c}
      N+2 \\ n+3 
    \end{array}\right) \frac{2^n}{n!} \int_0^\infty \rme^{-N\zeta}
  \zeta^n\,\rmd\zeta \\
&=  \sum_{n=0}^{N-1} \left(
    \begin{array}{c}
      N+2 \\ n+3 
    \end{array}\right) \frac{2^n}{N^{n+1}}.
\end{align*}
So,
\begin{equation}
  \tilde{F}_N({\textstyle\frac12},1) = \frac2\pi   \sum_{n=0}^{N-1} \left(
    \begin{array}{c}
      N+2 \\ n+3 
    \end{array}\right) \frac{2^n}{N^{n+1}},
\end{equation}
and
\begin{equation}
\tilde{F}({\textstyle\frac12},1) = \lim_{N\to\infty}
\frac1{N^2}  \tilde{F}_N({\textstyle\frac12},1)
=\lim_{N\to\infty}   \sum_{n=0}^{N-1} \left(
    \begin{array}{c}
      N+2 \\ n+3 
    \end{array}\right) \frac{2^n}{N^{n+3}}.
\end{equation}
To perform the last sum, we use the estimate $\displaystyle
\left(\begin{array}{c} N+2 \\ n+3 \end{array}\right) \leq
\frac{(N+2)^{n+3}}{(n+3)!}$ to get
\begin{align*}
  \left(\begin{array}{c} N+2 \\ n+3 \end{array}\right) 
  \frac{2^n}{N^{n+3}} & \leq
\frac{2^n}{(n+3)!}\left(\frac{N+2}N\right)^{n+3} \\
&\leq \frac{2^n}{n!}\rme^2,
\end{align*}
which is summable. This allows us to use Tannery's
theorem \cite[\S 49]{bro:ait}. Since
\begin{displaymath}
  \lim_{N\to\infty}\frac1{N^{n+3}}
    \left(\begin{array}{c} N+2 \\ n+3 \end{array}\right) 
 =\frac1{(n+3)!},
\end{displaymath}
by Stirling's approximation, we get
\begin{align}
  \nonumber \tilde{F}({\textstyle\frac12},1) &= 
  \frac2\pi \sum_{n=0}^\infty \frac{2^n}{(n+3)!} \\
  &=\frac{\rme^2-5}{4\pi}.
\end{align}
\section{Some sums}
In this appendix we give evaluations of some finite sums that
appear in the main text.
\subsection{A sum from section \ref{sec:some-integr-involv}}
\label{app:b_1}
In this appendix we evaluate the sum
\begin{equation}
  \label{eq:73}
  \sum_{\ell=0}^n (-1)^\ell \binomial{p}{n-\ell}
  \binomial{p-n+\ell-1}\ell,
\end{equation}
for $p\geq n+1$.

By elementary manipulation of binomial coefficients,
\begin{align}
  \label{eq:74}
  \binomial{p}{n-\ell}\binomial{p-n+\ell+1}\ell &=
  \frac{p!(p-n+\ell-1)!}{(n-\ell)!(p-n+\ell)!\ell!(p-n-1)!}
  \nonumber \\
  &=\frac{p}{p-n+\ell}\binomial{p-1}n\binomial{n}\ell.
\end{align}
We also use the fact that
\begin{equation}
  \label{eq:75}
  \int_0^1 x^{p-1-n+\ell}\,\rmd x = \frac1{p-n+\ell},
\end{equation}
together with \eqref{eq:74} to get
\begin{align}
  \label{eq:76}
  \sum_{\ell=0}^n (-1)^\ell \binomial{p}{n-\ell}
  \binomial{p-n+\ell-1}\ell &= p\binomial{p-1}n
  \sum_{\ell=0}^n (-1)^\ell \binomial{n}\ell
  \int_0^1 x^{p-1-n+\ell}\,\rmd x \nonumber \\
  &= p\binomial{p-1}n \int_0^1 x^{p-n-1}(1-x)^n\,\rmd x \nonumber \\
  &= p\binomial{p-1}n \frac{(p-n-1)!n!}{p!} = 1.
\end{align}
We have used
\begin{equation}
  \label{eq:77}
  \int_0^1 x^a(1-x)^b\,\rmd x = \frac{a!b!}{(a+b+1)!},\qquad\mbox{for
    $a, b\in\N$:}
\end{equation}
a version of Euler's integral.

\subsection{A sum from section \ref{sec:partition-sums-proof}}
\label{app:b_2}
In this appendix we evaluate the sum
\begin{equation}
  \label{eq:57}
  \sum_{x=0}^{p+1} \frac{(p-2x+1)^2}{x!(x+2)!
    (p-x+3)!(p-x+1)!}.
\end{equation}
Let us observe that the numerator in \eqref{eq:57} can be written
as
\begin{equation}
  \label{eq:58}
  (p-2x+1)^2 = (p-x+1)(p-x+3) + x(x+2) - (p-x+1)(x+2) - x(p-x+3).
\end{equation}
Therefore,
\begin{align}
  \label{eq:59}
  &\sum_{x=0}^{p+1} \frac{(p-2x+1)^2}{x!(x+2)!
    (p-x+3)!(p-x+1)!} \\ 
&= \sum_{x=0}^{p} \frac1{x!(x+2)!(p-x+2)!(p-x)!}
  + \sum_{x=1}^{p+1} \frac1{(x-1)!(x+1)!(p-x+3)!(p-x+1)!} \nonumber \\
&\qquad - \sum_{x=0}^{p} \frac1{x!(x+1)!(p-x+3)!(p-x)!}
- \sum_{x=1}^{p+1} \frac1{(x-1)!(x+2)!(p-x+2)!(p-x+1)!} \nonumber \\
&= \frac{1}{(p+2)!^2}\left( \sum_{x=0}^{p} \binomial{p+2}{x+2}
  \binomial{p+2}{x} + \sum_{x=1}^{p+1} \binomial{p+2}{x-1}
  \binomial{p+2}{x+1} \right) \nonumber \\
&\qquad -\frac{1}{(p+1)!(p+3)!} \left( \sum_{x=0}^{p}
  \binomial{p+1}{x+1}\binomial{p+3}x + \sum_{x=1}^{p+1}
  \binomial{p+1}{x-1}\binomial{p+3}{x+2}\right) \nonumber \\
&= \frac{2}{(p+2)!^2} \sum_{x=0}^{p} \binomial{p+2}{x}
  \binomial{p+2}{p-x}  \nonumber \\
&\qquad -\frac{1}{(p+1)!(p+3)!} \left( \sum_{x=0}^{p}
  \binomial{p+3}{x}\binomial{p+1}{p-x} + \sum_{x=0}^{p}
  \binomial{p+1}{x}\binomial{p+3}{p-x}\right). \nonumber 
\end{align}
We then use the Vandermonde identity,
\begin{equation}
  \label{eq:60}
  \sum_{x=0}^{p} \binomial{w}{x}\binomial{v}{p-x} = 
  \binomial{w+v}{p},
\end{equation}
to get
\begin{align}
  \label{eq:61}
  \sum_{x=0}^{p+1} \frac{(p-2x+1)^2}{x!(x+2)!
    (p-x+3)!(p-x+1)!} &= \frac2{(p+2)!^2}\binomial{2p+4}{p} \nonumber\\
 &\qquad -\frac2{(p+1)!(p+3)!}\binomial{2p+4}p \nonumber \\
  &=\frac{2}{(p+1)!(p+2)!}\left(\frac1{p+2}-\frac1{p+3}\right)
  \binomial{2p+4}p \nonumber\\
  &=\frac2{(p+2)!(p+3)!}\binomial{2p+4}p.
\end{align}
\def\Dbar{\leavevmode\lower.6ex\hbox to 0pt{\hskip-.23ex \accent"16\hss}D}
  \def\cprime{$'$} \def\rmi{{\mathrm i}}

\end{document}